\documentclass[reqno]{amsart}
\usepackage{hyperref}

\AtBeginDocument{{\noindent\small
\emph{}}
\vspace{8mm}}
\begin{document}
\title[\hfil fractional Schr\"{o}dinger-Poisson system]
{Concentrating bounded states for fractional Schr\"{o}dinger-Poisson system involving critical Sobolev exponent}

\author[K. M. Teng ]
{Kaimin Teng}  
\address{Kaimin Teng (Corresponding Author)\newline
Department of Mathematics, Taiyuan
University of Technology, Taiyuan, Shanxi 030024, P. R. China}
\email{tengkaimin2013@163.com}


\begin{abstract}
In this paper, we study the concentration and multiplicity of solutions to the following fractional Schr\"{o}dinger-Poisson system
\begin{equation*}
\left\{
  \begin{array}{ll}
    \varepsilon^{2s}(-\Delta)^su+V(x)u+\phi u=f(u)+u^{2_s^{\ast}-1} & \hbox{in $\mathbb{R}^3$,} \\
     \varepsilon^{2t}(-\Delta)^t\phi=u^2, u>0& \hbox{in $\mathbb{R}^3$,}
  \end{array}
\right.
\end{equation*}
where $s>\frac{3}{4}$, $s,t\in(0,1)$, $\varepsilon>0$ is a small parameter, $f\in C^1(\mathbb{R}^{+},\mathbb{R})$ is subcritical, $V:\mathbb{R}^3\rightarrow\mathbb{R}$ is a continuous bounded function. We establish a family of positive solutions $u_{\varepsilon}\in H_{\varepsilon}$ which concentrates around the local minima of $V$ in $\Lambda$ as $\varepsilon\rightarrow0$. With Ljusternik-Schnirelmann theory, we also obtain multiple solutions by employing the topology construct of the set where the potential $V$ attains its minimum.
\end{abstract}

\maketitle
\numberwithin{equation}{section}
\newtheorem{theorem}{Theorem}[section]
\newtheorem{lemma}[theorem]{Lemma}
\newtheorem{definition}[theorem]{Definition}
\newtheorem{remark}[theorem]{Remark}
\newtheorem{proposition}[theorem]{Proposition}
\newtheorem{corollary}[theorem]{Corollary}
\allowdisplaybreaks

\section{Introduction}
In this paper, we study the following fractional Schr\"{o}dinger-Poisson system
\begin{equation}\label{main}
\left\{
  \begin{array}{ll}
    \varepsilon^{2s}(-\Delta)^su+V(x)u+\phi u=f(u)+u^{2_s^{\ast}-1} & \hbox{in $\mathbb{R}^3$,} \\
     \varepsilon^{2t}(-\Delta)^t\phi=u^2, u>0& \hbox{in $\mathbb{R}^3$,}
  \end{array}
\right.
\end{equation}
where $s>\frac{3}{4}$, $s,t\in(0,1)$, $\varepsilon>0$ is a small parameter. The potential $V:\mathbb{R}^3\rightarrow\mathbb{R}$ is a continuous bounded function satisfying\\
$(V_0)$ $\inf\limits_{x\in\mathbb{R}^3}V(x)=\alpha_0>0$;\\
$(V_1)$ There is a bounded domain $\Lambda\subset\mathbb{R}^3$ such that
\begin{equation*}
V_0=\inf_{\Lambda}V(x)<\min_{\partial\Lambda}V(x).
\end{equation*}
This kind of hypothesis was first introduced by del Pino and Felmer in \cite{DF}. The nonlinearity $f:\mathbb{R}\rightarrow\mathbb{R}$ is of $C^1$-class function. Since we are looking for positive solutions, we may assume that $f(\tau)=0$ for $\tau<0$. Furthermore, we need the following conditions:\\
$(f_0)$ $\lim\limits_{\tau\rightarrow0^{+}}\frac{f(\tau)}{\tau^3}=0$;\\
$(f_1)$ there exist $\lambda>0$ and $C>0$ such that $f(\tau)\geq\lambda \tau^{q-1}$ for some $4\leq q<2_s^{\ast}$ and $|f'(\tau)|\leq C(1+|\tau|^{p-2})$, where $4<p<2_s^{\ast}$ ;\\
$(f_2)$ $\frac{f(\tau)}{\tau^{3}}$ is non-decreasing in $\tau\in (0,+\infty)$.

The non-local operator $(-\Delta)^s$ ($s\in(0,1)$), which is called fractional Laplacian operator, can be defined by
\begin{equation*}
(-\Delta)^su(x)=C_s\,{\rm P.V.}\int_{\mathbb{R}^3}\frac{u(x)-u(y)}{|x-y|^{3+2s}}\,{\rm d}y=C_s\lim_{\varepsilon\rightarrow0}\int_{\mathbb{R}^3\backslash B_{\varepsilon}(x)}\frac{u(x)-u(y)}{|x-y|^{3+2s}}\,{\rm d}y
\end{equation*}
for $u\in\mathcal{S}(\mathbb{R}^3)$, where $\mathcal{S}(\mathbb{R}^3)$ is the Schwartz space of rapidly decaying $C^{\infty}$ function, $B_{\varepsilon}(x)$ denote an open ball of radius $r$ centered at $x$ and the normalization constant $C_s=\Big(\int_{\mathbb{R}^3}\frac{1-\cos(\zeta_1)}{|\zeta|^{3+2s}}\,{\rm d}\zeta\Big)^{-1}$. For $u\in\mathcal{S}(\mathbb{R}^3)$, the fractional Laplace operator $(-\Delta)^s$ can be expressed as an inverse Fourier transform
\begin{equation*}
(-\Delta)^su=\mathcal{F}^{-1}\Big((2\pi|\xi|)^{2s}\mathcal{F}u(\xi)\Big),
\end{equation*}
where $\mathcal{F}$ and $\mathcal{F}^{-1}$ denote the Fourier transform and inverse transform, respectively. If $u$ is sufficiently smooth, it is known that (see \cite{NPV}) it is equivalent to
\begin{equation*}
(-\Delta)^su(x)=-\frac{1}{2}C_s\int_{\mathbb{R}^3}\frac{u(x+y)+u(x-y)-2u(y)}{|x-y|^{3+2s}}\,{\rm d}y.
\end{equation*}
By a classical solution of \eqref{main}, we mean two continuous functions that $(-\Delta)^su$ is well defined for all $x\in\mathbb{R}^3$ and satisfies \eqref{main} in pointwise sense.

In the last several years, nonlinear equations involving fractional Laplacian is much of interest, and attracts much attention by many scholars. One of the main reason is that fractional operators appear in many mathematical and physical problems, such as: fractional quantum mechanics \cite{La,La1}, Financial modelling \cite{CT}, anomalous diffusion \cite{MK}, obstacle problems \cite{S}, conformal geometry and minimal surfaces \cite{CM} and so on. Another main reason is that the fractional Laplacian $(-\Delta)^s$ ($s\in(0,1)$) is a non-local operator comparing with the classical Laplacian $-\Delta$ which is a local one, the methods which were previously developed, maybe not be applied directly. We refer the interesting readers to see the recent progresses such as \cite{BCPS,CS,CW,DPW,FQT,FL,MMS,MS,NPV,S,Sec,SZ,Teng,Teng1,Teng2,ZDS} and the references therein.

In the very recent, fractional Schr\"{o}dinger-Poisson system
\begin{equation}\label{main1}
\left\{
  \begin{array}{ll}
   (-\Delta)^su+V(x)u+\phi u=f(u) & \hbox{in $\mathbb{R}^3$,} \\
   (-\Delta)^t\phi=u^2& \hbox{in $\mathbb{R}^3$,}
  \end{array}
\right.
\end{equation}
has been investigated by some scholars. When $f(u)=|u|^{p-1}u$ for $2<p<2_s^{\ast}-1$ or $f(u)=\mu|u|^{q-1}+|u|^{2_s^{\ast}-2}u$ with $q\in(\frac{3s+t}{s+t},2_s^{\ast})$ and $\mu>0$ which maybe large for some $q$, in \cite{Teng,Teng1}, we established the existence of positive ground state solution by using the Nehari-Pohozaev manifold combing monotone trick with global compactness Lemma, respectively. In \cite{ZDS}, the authors studied the existence of radial solutions for system \eqref{main1} with the nonlinearity $f(u)$ verifying the subcritical or critical assumptions of Berestycki-Lions type.

For the semiclassical state, in \cite{MS}, the authors studied the semiclassical state of the following system
\begin{equation*}
\left\{
  \begin{array}{ll}
    \varepsilon^{2s}(-\Delta)^su+V(x)u+\phi u=f(u) & \hbox{in $\mathbb{R}^N$,} \\
    \varepsilon^{\theta}(-\Delta)^{\frac{\alpha}{2}}\phi=\gamma_{\alpha}u^2& \hbox{in $\mathbb{R}^N$,}
  \end{array}
\right.
\end{equation*}
where $s\in(0,1)$, $\alpha\in(0,N)$, $\theta\in(0,\alpha)$, $N\in(2s,2s+\alpha)$, $\gamma_{\alpha}$ is a positive constant, $f(u)$ satisfies the following subcritical growth assumptions: $0<KF(t)\leq f(t)t$ with some $K>4$ for all $t\geq0$ and $\frac{f(t)}{t^3}$ is strictly increasing on $(0,+\infty)$. By adapting some ideas of Benci, Cerami and Passaseo \cite{BC,BCP} and using the Ljusternick-Schnirelmann Theory, the authors obtained the multiplicity of positive solutions which concentrate on the minima of $V(x)$ as $\varepsilon\rightarrow0$. Using the similar methods, the authors in \cite{LZ} studied the system \eqref{main} and established the multiplicity and concentration behavior of solutions. In \cite{Teng2}, we also studied the concentration of positive ground state solution via the Nehari manifold for the following system
\begin{equation*}
\left\{
  \begin{array}{ll}
    \varepsilon^{2s}(-\Delta)^su+V(x)u+\phi u=K(x)f(u)+Q(x)|u|^{2_s^{\ast}-2}u & \hbox{in $\mathbb{R}^3$,} \\
     \varepsilon^{2t}(-\Delta)^t\phi=u^2& \hbox{in $\mathbb{R}^3$,}
  \end{array}
\right.
\end{equation*}
where $V(x)$, $K(x)$, $Q(x)$ is a bounded continuous potential and $f$ satisfies the subcritical growth and some monotone condition. In the above works, the potential $V(x)$ is either constant or some bounded potentials possessing some global minimum points, the main purpose of this paper is devoted to studying the case that $V(x)$ possesses some local minimum points.

In the last several years, the semiclassical state of the nonlinear Schr\"{o}dinger-Poisson system has been object of interest for many authors. Ruiz and Vaira \cite{RV} proved the existence of multi-bump solutions of system
\begin{equation}\label{main-1-1}
\left\{
  \begin{array}{ll}
    -\varepsilon^2\Delta u+V(x)u+K(x)\phi u=u^p & \hbox{in $\mathbb{R}^3$,} \\
    -\Delta \phi=K(x)u^2& \hbox{in $\mathbb{R}^3$.}
  \end{array}
\right.
\end{equation}
with $K(x)\equiv1$ and these bumps concentrate around a local minimum of the potential $V$. Ruiz \cite{Ruiz} and D'Aprile and Wei \cite{DW} showed that system \eqref{main-1-1} with
$V(x)\equiv K(x)\equiv1$ possesses a family of solutions concentrating around a sphere when $\varepsilon\rightarrow0$ for $p\in(2,\frac{18}{7})$. Their results were generalized in \cite{IV1,I} for the radial $V$ and $K$. Ianni and Vaira \cite{IV} obtained the existence of positive bound state solutions which concentrate on a non-degenerate local minimum or maximum of $V$ by using a Lyapunov-Schmitt reduction method. Seok \cite{Seok} proved that system \eqref{main-1-1} has single and multi-peak solutions which concentrate around the local minimum of $V$ with the Berestycki-Lions conditions. After that, Zhang \cite{Zhang} considered the critical Berestycki-Lions conditions and obtained the single solutions which concentrate around the local minimum of $V$. Liu at el \cite{LGF} proved the multi-semiclassical states which concentrates around its corresponding global minimum point of $V$.
The concentration phenomenon of solutions for the following Schr\"{o}dinger-Poisson system also investigated by some authors:
\begin{equation}\label{main-1-2}
\left\{
  \begin{array}{ll}
    -\varepsilon^2\Delta u+V(x)u+\lambda\phi u=f(x,u) & \hbox{in $\mathbb{R}^3$,} \\
    -\varepsilon^2\Delta \phi=u^2& \hbox{in $\mathbb{R}^3$.}
  \end{array}
\right.
\end{equation}
For the subcritical case, When $f(x,u)=b(x)f(u)$ and $f$ satisfies super-4 growth condition and some monotone condition, Wang et al. \cite{WTXZ} studied the existence of ground state solutions for system \eqref{main-1-2} and the concentration behavior of least energy solutions was obtained. For the critical growth case, He and Zou \cite{HZ} considered the existence and concentration behavior of ground state solutions for \eqref{main-1-2} with $f(x,u)=u^5+ f(u)$ and the subcritical term $f$ the so-called Ambrosetti-Rabinowitz condition, they proved that system \eqref{main-1-2} has a ground state solution concentrating around a global minimum $V$ as $\varepsilon\rightarrow0$. When $f(x,u)=\lambda|u|^{p-2}u+|u|^4u$ with $3<p\leq 4$, which it does not satisfy the monotone assumption or Ambrosetti-Rabinowtiz condition, He and Li \cite{HL} construct a family of positive solution which concentrates around a local minimum of $V$ as $\varepsilon\rightarrow0$.

But, to the best of knowledge, for the local case, only in \cite{HL,Seok, Zhang}, the authors considered the the potential $V(x)$ possessing a local minimum point. They followed the method developed by Byeon-Jeanjean \cite{BJ,BJT} to construct a peak or multiple-peak solution which concentrates on the local minimum of $V(x)$ as $\varepsilon\rightarrow0$. For the nonlocal case, there are no papers to consider the potential $V(x)$ possessing a local minimum point. Motivated by some related works, the aim of this paper is to study the existence of concentration solutions in the case that $V(x)$ has local minimum points. We will take the penalization arguments due to del Pino and Felmer \cite{PF} to investigate system \eqref{main}. As we know, this kind of penalization method has been successfully applied to study the multiplicity and concentration of solutions for other problems, such as: Kirchhoff type problems \cite{HLP}, fractional Schr\"{o}dinger equations \cite{AM}, quasilinear problem involving $N$-Laplacian \cite{ASA}, quasilinear Choquard equation \cite{AY} and so on. However, now we are working with a class of non-local problems, there are some new difficulties in dealing with the system \eqref{main}. One of the main difficulties is that we consider the system \eqref{main} being with critical Sobolev exponent, it is required to use the concentration-compactness principle to return the compactness, but in \cite{PP}, the authors provided a version of concentration-compactness principle which is useful for the bounded domain, and here our problem is set on the whole space $\mathbb{R}^3$. This is the first obstacle to be killed. Another is the decay estimate of solution sequence at infinity, this is different from the classical one, such as \eqref{main-1-2}. These difficulties make us more careful analysis, which permit us to use the penalization method.

Our main results are as follows.
\begin{theorem}\label{thm1-1}
Let $2s+2t>3$, $s,t\in(0,1)$ and $s>\frac{3}{4}$. Suppose that $V$ satisfies $(V_0)$, $(V_1)$ and $f\in C(\mathbb{R}^{+},\mathbb{R})$ satisfies $(f_0)$-$(f_2)$. Then there exists an $\varepsilon_0>0$ such that system \eqref{main} possesses a positive solution $(u_{\varepsilon},\phi_{\varepsilon})\in H_{\varepsilon}\times\mathcal{D}^{t,2}(\mathbb{R}^3)$ for all $\varepsilon\in(0,\varepsilon_0)$. Moreover, $u_{\varepsilon}$ possesses a maximum $x_{\varepsilon}\in \Lambda$ such that $V(x_{\varepsilon})\rightarrow \inf\limits_{\Lambda}V$ as $\varepsilon\rightarrow0$, and
\begin{equation*}
u_{\varepsilon}(x)\leq\frac{C\varepsilon^{3+2s}}{C_0\varepsilon^{3+2s}+|x-x_{\varepsilon}|^{3+2s}}\quad x\in\mathbb{R}^3,\,\,\text{and}\,\, \varepsilon\in(0,\varepsilon_0)
\end{equation*}
for some constants $C>0$ and $C_0\in\mathbb{R}$.
\end{theorem}
By using the Lyusternik-Schnirelmann category, we can obtain the multiplicity of positive solutions. For this purpose, we also need the assumption that
\begin{equation}
(V_2)\quad \mathcal{M}=\{x\in\Lambda\,\,|\,\, V(x)=\inf_{\mathbb{R}^3}V\}\neq\emptyset.
\end{equation}
We denote the closed $\delta$-neighborhood of the set $\mathcal{M}$ by
\begin{equation*}
\mathcal{M}_{\delta}=\{x\in\mathbb{R}^3\,\,|\,\,{\rm dist}(x,\mathcal{M})\leq\delta\}.
\end{equation*}
Also we recall that, if $Y$ is a closed set of a topological space $X$, ${\rm cat}_{X}(Y)$ is the Ljusternik-Schnirelmann category of $Y$ in $X$, namely the least number of closed and contractible sets in $X$ which cover $Y$. We shall prove the following multiplicity result.
\begin{theorem}\label{thm1-2}
Let $2s+2t>3$, $s,t\in(0,1)$ and $s>\frac{3}{4}$. Suppose that $V$ satisfies $(V_0)$, $(V_1)$, $(V_2)$ and $f\in C(\mathbb{R}^{+},\mathbb{R})$ satisfies $(f_0)$-$(f_2)$. Then, for any $\delta>0$ given, there exists $\varepsilon_{\delta}>0$ such that, for any $\varepsilon\in(0,\varepsilon_{\delta})$, system \eqref{main} has at least ${\rm cat}_{\mathcal{M}_{\delta}}(\mathcal{M})$ solutions. Furthermore, if one of these solutions $u_{\varepsilon}$ possesses a maximum $x_{\varepsilon}\in\Lambda$, then $\lim\limits_{\varepsilon\rightarrow0}V(x_{\varepsilon})=\inf\limits_{\Lambda}V$ and
\begin{equation*}
u_{\varepsilon}(x)\leq\frac{C\varepsilon^{3+2s}}{C_0\varepsilon^{3+2s}+|x-x_{\varepsilon}|^{3+2s}}\quad x\in\mathbb{R}^3,\,\,\text{and}\,\, \varepsilon\in(0,\varepsilon_{\delta})
\end{equation*}
for some constants $C>0$ and $C_0\in\mathbb{R}$.
\end{theorem}

This paper is organized as follows, in section 2, we give some preliminary results and a version of concentration-compactness principle. In section 3, we will prove the penalization problem has a positive solution. In section 4, we will prove the limiting problem has a positive ground state solution. In section 5, we will give a uniform estimate for solution sequences. In section 6, we complete the proof of Theorem \ref{thm1-1}. Section 7 is devote to prove Theorem \ref{thm1-2} by using the Lyusternik-Schnirelmann category.

\section{Variational Setting}

In this section, we outline the variational framework for studying problem \eqref{main} and list some preliminary Lemma which used later. In the sequel, we denote by $\|\cdot\|_{p}$ the usual norm of the space $L^p(\mathbb{R}^3)$, the letter $c_i$ ($i=1,2,\ldots$) or $C$ denote by some positive constants. We denote $\widehat{u}$ the Fourier transform of $u$.

\subsection{Work space stuff}
We define the homogeneous fractional Sobolev space $\mathcal{D}^{s,2}(\mathbb{R}^3)$ as follows
\begin{equation*}
\mathcal{D}^{s,2}(\mathbb{R}^3)=\Big\{u\in L^{2_{s}^{\ast}}(\mathbb{R}^3)\,\,\Big|\,\,|\xi|^{s}\widehat{u}(\xi)\in L^2(\mathbb{R}^3)\Big\}
\end{equation*}
which is the completion of $C_0^{\infty}(\mathbb{R}^3)$ under the norm
\begin{equation*}
\|u\|_{\mathcal{D}^{s,2}}=\Big(\int_{\mathbb{R}^3}|(-\Delta)^{\frac{s}{2}}u|^2\,{\rm d}x\Big)^{\frac{1}{2}}=\Big(\int_{\mathbb{R}^3}|\xi|^{2s}|\widehat{u}(\xi)|^2\,{\rm d}\xi\Big)^{\frac{1}{2}}
\end{equation*}
The fractional Sobolev space $H^{s}(\mathbb{R}^3)$ can be described by means of the Fourier transform, i.e.
\begin{equation*}
H^{s}(\mathbb{R}^3)=\Big\{u\in L^2(\mathbb{R}^3)\,\,\Big|\,\,\int_{\mathbb{R}^3}(|\xi|^{2s}|\widehat{u}(\xi)|^2+|\widehat{u}(\xi)|^2)\,{\rm d}\xi<+\infty\Big\}.
\end{equation*}
In this case, the inner product and the norm are defined as
\begin{equation*}
(u,v)=\int_{\mathbb{R}^3}(|\xi|^{2s}\widehat{u}(\xi)\overline{\widehat{v}(\xi)}+\widehat{u}(\xi)\overline{\widehat{v}(\xi)})\,{\rm d}\xi
\end{equation*}
\begin{equation*}
\|u\|_{H^{s}}=\bigg(\int_{\mathbb{R}^3}(|\xi|^{2s}|\widehat{u}(\xi)|^2+|\widehat{u}(\xi)|^2)\,{\rm d}\xi\bigg)^{\frac{1}{2}},
\end{equation*}
From Plancherel's theorem we have $\|u\|_2=\|\widehat{u}\|_2$ and $\||\xi|^{s}\widehat{u}\|_2=\|(-\Delta)^{\frac{s}{2}}u\|_2$. Hence
\begin{equation*}
\|u\|_{H^{s}}=\bigg(\int_{\mathbb{R}^3}(|(-\Delta)^{\frac{s}{2}}u(x)|^2+|u(x)|^2)\,{\rm d}x\bigg)^{\frac{1}{2}},\quad \forall u\in H^{s}(\mathbb{R}^3).
\end{equation*}
We denote $\|\cdot\|$ by $\|\cdot\|_{H^{s}}$ in the sequel for convenience.

We define the Sobolev space $H_{\varepsilon}=\{u\in H^s(\mathbb{R}^3)\,\,|\,\, \int_{\mathbb{R}^3}V(\varepsilon x)u^2\,{\rm d}x<\infty\}$ endowed with the norm
\begin{equation*}
\|u\|_{H_{\varepsilon}}=\Big(\int_{\mathbb{R}^3}\Big(|(-\Delta)^{\frac{s}{2}}u|^2+V(\varepsilon x)u^2\Big)\,{\rm d}x\Big)^{\frac{1}{2}}.
\end{equation*}

It is well known that $H^{s}(\mathbb{R}^3)$ is continuously embedded into $L^r(\mathbb{R}^3)$ for $2\leq r\leq 2_{\alpha}^{\ast}$ ($2_{\alpha}^{\ast}=\frac{6}{3-2\alpha}$).
Obviously, the conclusion also holds for $H_{\varepsilon}$.

\subsection{Formulation of Problem \eqref{main}}
It is easily seen that, just performing the change of variables $u(x)\rightarrow u(x/\varepsilon)$ and $\phi(x)\rightarrow \phi(x/\varepsilon)$, and taking $z=x/\varepsilon$, problem \eqref{main} can be rewritten as the following equivalent form
\begin{equation}\label{main1}
\left\{
  \begin{array}{ll}
   (-\Delta)^su+V(\varepsilon z)u+\phi u=f(u)+u^{2_s^{\ast}-1} & \hbox{in $\mathbb{R}^3$,} \\
    (-\Delta)^t\phi=u^2, u>0& \hbox{in $\mathbb{R}^3$}
  \end{array}
\right.
\end{equation}
which will be referred from now on.

From $(f_2)$, it is easy to verify that
\begin{equation*}
f'(\tau)\tau-3f(\tau)>0\,\,\text{and}\,\,  f(\tau)\tau-4F(\tau)>0 \,\,\text{for any}\,\,\tau>0.
\end{equation*}

Observe that if $4s+2t\geq3$, there holds $2\leq\frac{12}{3+2t}\leq\frac{6}{3-2s}$ and thus $H_{\varepsilon}\hookrightarrow L^{\frac{12}{3+2t}}(\mathbb{R}^3)$. Considering $u\in H_{\varepsilon}$, the linear functional $\widetilde{\mathcal{L}}_u:\mathcal{D}^{t,2}(\mathbb{R}^3)\rightarrow\mathbb{R}$ is defined by
\begin{equation*}
\widetilde{\mathcal{L}}_u(v)=\int_{\mathbb{R}^3}u^2v\,{\rm d}z.
\end{equation*}
Using the Lax-Milgram theorem, there exists a unique $\phi_u^t\in\mathcal{D}^{t,2}(\mathbb{R}^3)$ such that
\begin{equation*}
\int_{\mathbb{R}^3}(-\Delta)^{\frac{t}{2}}\phi_u^t(-\Delta)^{\frac{t}{2}}v\,{\rm d}z=\int_{\mathbb{R}^3}u^2v\,{\rm d}z,\quad \forall v\in\mathcal{D}^{t,2}(\mathbb{R}^3),
\end{equation*}
that is $\phi_u^t$ is a weak solution of $(-\Delta)^t\phi_u^t=u^2$ and so the representation formula holds
\begin{equation*}
\phi_u^t(z)=c_t\int_{\mathbb{R}^3}\frac{u^2(y)}{|z-y|^{3-2t}}\,{\rm d}y,\quad z\in\mathbb{R}^3,\quad c_t=\pi^{-\frac{3}{2}}2^{-2t}\frac{\Gamma(\frac{3-2t}{2})}{\Gamma(t)}.
\end{equation*}
Substituting $\phi_u^t$ in \eqref{main1}, it reduces to a single fractional Schr\"{o}dinger equation
\begin{equation}\label{R-1}
(-\Delta)^su+V(\varepsilon z)u+\phi_u^tu=f(u)+u^{2_s^{\ast}-1},\quad z\in\mathbb{R}^3.
\end{equation}
The solvation of \eqref{R-1} can be looking for the critical points of the associated energy functional $J_{\varepsilon}: H_{\varepsilon}\rightarrow\mathbb{R}$ defined by
\begin{equation*}
J_{\varepsilon}(u)=\frac{1}{2}\int_{\mathbb{R}^3}\Big(|(-\Delta)^{\frac{s}{2}}u|^2+V(\varepsilon z)u^2\Big)\,{\rm d}z+\frac{1}{4}\int_{\mathbb{R}^3}\phi_u^tu^2\,{\rm d}z-\int_{\mathbb{R}^3}F(u)\,{\rm d}z-\int_{\mathbb{R}^3}(u^{+})^{2_s^{\ast}}\,{\rm d}z
\end{equation*}
and $J_{\varepsilon}\in C^1(H_{\varepsilon},\mathbb{R})$.
Let us summarize some properties of the function $\phi_u^t$, the proof can be found in \cite{Teng2}.
\begin{lemma}\label{lem2-1}
For every $u\in H_{\varepsilon}$ with $4s+2t\geq3$, define $\Phi(u)=\phi_u^t\in \mathcal{D}^{t,2}(\mathbb{R}^3)$, where $\phi_u^t$ is the unique solution of equation $(-\Delta)^t\phi=u^2$. Then there hold:\\
$(i)$ If $u_n\rightharpoonup u$ in $H_{\varepsilon}$, then $\Phi(u_n)\rightharpoonup\Phi(u)$ in $\mathcal{D}^{t,2}(\mathbb{R}^3)$;\\
$(ii)$ $\Phi(tu)=t^2\Phi(u)$ for any $t\in\mathbb{R}$;\\
$(iii)$ For $u\in H_{\varepsilon}$, one has
\begin{equation*}
\|\Phi(u)\|_{\mathcal{D}^{t,2}}\leq C\|u\|_{\frac{12}{3+2t}}^2\leq C\|u\|_{H_\varepsilon}^2,\quad \int_{\mathbb{R}^3}\Phi(u)u^2\,{\rm d}x\leq C\|u\|_{\frac{12}{3+2t}}^4\leq C\|u\|_{H_\varepsilon}^4,
\end{equation*}
where constant $C$ is independent of $u$;\\
$(iv)$ Let $2s+2t>3$, if $u_n\rightharpoonup u$ in $H_{\varepsilon}$ and $u_n\rightarrow u$ a.e. in $\mathbb{R}^3$, then for any $v\in H_{\varepsilon}$,
\begin{equation*}
\int_{\mathbb{R}^3}\phi_{u_n}^tu_nv\,{\rm d}z\rightarrow\int_{\mathbb{R}^3}\phi_{u}^tuv\,{\rm d}z\quad\text{and}\quad\int_{\mathbb{R}^3}f(u_n)v\,{\rm d}z\rightarrow\int_{\mathbb{R}^3}f(u)v\,{\rm d}z
\end{equation*}
and
\begin{equation*}
\int_{\mathbb{R}^3}(u_n^{+})^{2_s^{\ast}-1}v\,{\rm d}z\rightarrow\int_{\mathbb{R}^3}(u^{+})^{2_s^{\ast}-1}v\,{\rm d}z.
\end{equation*}
\end{lemma}

In the end of this section, we will give a version of concentration-compactness on whole space $\mathbb{R}^3$ which is sufficient to prove our main results.
We define
\begin{equation*}
\mu_{\infty}=\lim_{R\rightarrow\infty}\limsup_{n\rightarrow\infty}\int_{\mathbb{R}^3\backslash B_R(0)}|(-\Delta)^{\frac{s}{2}}u_n|^2\,{\rm d}z\quad \nu_{\infty}=\lim_{R\rightarrow\infty}\limsup_{n\rightarrow\infty}\int_{\mathbb{R}^3\backslash B_R(0)}|u_n|^{2_s^{\ast}}\,{\rm d}z.
\end{equation*}
\begin{lemma}\label{lem2-2}(Theorem 5 \cite{PP}, Lemma 3.5 and 3.7 \cite{ZZX})
Let $\{u_n\}\subset \mathcal{D}^{s,2}(\mathbb{R}^3)$ be such that $u_n\rightharpoonup u$ in $\mathcal{D}^{s,2}(\mathbb{R}^3)$, $u_n\rightarrow u$ in $L^2(\mathbb{R}^3)$, $|(-\Delta)^{\frac{s}{2}}u_n|^2\rightharpoonup\mu$ and $|u_n|^{2_s^{\ast}}\rightharpoonup\nu$ weakly$-\ast$ in $\mathcal{M}(\mathbb{R}^3)$ as $n\rightarrow\infty$. Here $\mathcal{M}(\mathbb{R}^3)$ is the space of finite nonnegative Borel measures on $\mathbb{R}^3$. Then\\
$(i)$ $u_n\rightarrow u$ in $L_{loc}^{2_s^{\ast}}$ or there exists a (at most countable) set of distinct points $\{x_j\}_{j\in J}\subset\mathbb{R}^3$ and positive number $\{\nu_j\}_{j\in J}$ such that
\begin{equation*}
\nu=|u|^{2_s^{\ast}}+\sum_{j\in J}\nu_j\delta_{x_j};
\end{equation*}
$(ii)$
Then $\mu_{\infty}$ and $\nu_{\infty}$ are well defined satisfy
\begin{equation*}
\limsup_{n\rightarrow\infty}\int_{\mathbb{R}^3}|(-\Delta)^{\frac{s}{2}}u_n|^2\,{\rm d}z=\int_{\mathbb{R}^3}\,{\rm d}\mu+\mu_{\infty}\quad \limsup_{n\rightarrow\infty}\int_{\mathbb{R}^3}|u_n|^{2_s^{\ast}}\,{\rm d}z=\int_{\mathbb{R}^3}\,{\rm d}\nu+\nu_{\infty};
\end{equation*}
$(iii)$
\begin{equation*}
\nu_j\leq(\mathcal{S}_s^{-1}\mu_j(\{x_j\}))^{\frac{2_s^{\ast}}{2}}\,\, \text{for any}\,\, j\in J\,\, \text{and}\,\,\nu_{\infty}\leq(\mathcal{S}_s^{-1}\mu_{\infty})^{\frac{2_s^{\ast}}{2}}.
\end{equation*}
\end{lemma}
\begin{proof}
The conclusion $(i)$ comes from Theorem 5 in \cite{PP}, $(ii)$ comes from Lemma 3.5 in \cite{ZZX}. We only need to show that $(iii)$ holds.

1. Take $\psi\in C_0^{\infty}(\mathbb{R}^3)$ such that $\psi=1$ on $B_1(0)$, $\psi=0$ on $\mathbb{R}^3\backslash B_2(0)$, $0\leq|\psi\leq1$ and $|\nabla\psi|\leq C$. For any $\rho>0$, define $\psi_{\rho}(x)=\psi(\frac{x-x_j}{\rho})$, where $j\in J$. It follows from Sobolev inequality that
\begin{equation*}
\int_{\mathbb{R}^3}|u_n\psi_{\rho}|^{2_s^{\ast}}\,{\rm d}x\leq\Big(\mathcal{S}_s^{-1}\int_{\mathbb{R}^3}|(-\Delta)^{\frac{s}{2}}(u_n\psi_{\rho})|^2\,{\rm d}y\,{\rm d}x\Big)^{\frac{2_s^{\ast}}{2}}.
\end{equation*}
Since $|u_n|^{2_s^{\ast}}\rightharpoonup \nu$ in $\mathcal{M}(\mathbb{R}^3)$, we have
\begin{equation*}
\int_{\mathbb{R}^3}|u_n\psi_{\rho}|^{2_s^{\ast}}\,{\rm d}x\rightarrow\int_{\mathbb{R}^3}\psi_{\rho}^{2_s^{\ast}}\,{\rm d}\nu\rightarrow\nu(\{x_j\})=\nu_j\quad \text{as}\,\, \rho\rightarrow0.
\end{equation*}
Using the nonlocal Leibniz rule:
\begin{equation*}
(-\Delta)^{\frac{s}{2}}(u_n\psi_{\rho})=\psi_{\rho}(-\Delta)^{\frac{s}{2}}u_n+u_n(-\Delta)^{\frac{s}{2}}\psi_{\rho}-B(u_n,\psi_{\rho}),
\end{equation*}
where $B(u_n,\psi_{\rho})=C_s{\rm P.V.}\int_{\mathbb{R}^3}\frac{(\psi_{\rho}(x)-\psi_{\rho}(y))(u_n(x)-u_n(y))}{|x-y|^{3+s}}\,{\rm d}y$. Then it is easy to obtain
\begin{align*}
&\int_{\mathbb{R}^3}|(-\Delta)^{\frac{s}{2}}(u_n\psi_{\rho})|^2\,{\rm d}x=\int_{\mathbb{R}^3}|(-\Delta)^{\frac{s}{2}}u_n|^2\psi_{\rho}^2\,{\rm d}x+\int_{\mathbb{R}^3}u_n^2|(-\Delta)^{\frac{s}{2}}\psi_{\rho}|^2\,{\rm d}x+\int_{\mathbb{R}^3}B^2(u_n,\psi_{\rho})\,{\rm d}x\\
&+2\int_{\mathbb{R}^3}u_n\psi_{\rho}(-\Delta)^{\frac{s}{2}}u_n(-\Delta)^{\frac{s}{2}}\psi_{\rho}\,{\rm d}x-2\int_{\mathbb{R}^3}\psi_{\rho}(-\Delta)^{\frac{s}{2}}u_nB(u_n,\psi_{\rho})\,{\rm d}x\\
&-2\int_{\mathbb{R}^3}u_n(-\Delta)^{\frac{s}{2}}\psi_{\rho}B(u_n,\psi_{\rho})\,{\rm d}x:=\int_{\mathbb{R}^3}|(-\Delta)^{\frac{s}{2}}u_n|^2\psi_{\rho}^2\,{\rm d}x+A_1+A_2+A_3+A_4+A_5.
\end{align*}
Next, we will show that $\lim\limits_{\rho\rightarrow0^{+}}\limsup\limits_{n\rightarrow\infty}A_i=0$, $i=1,2,3,4,5$. If these are proved, using the assumption $|(-\Delta)^{\frac{s}{2}}u_n|^2\rightharpoonup\mu$ in $\mathcal{M}(\mathbb{R}^3)$, then
\begin{equation*}
\lim_{\rho\rightarrow0^{+}}\limsup_{n\rightarrow\infty}\int_{\mathbb{R}^3}|(-\Delta)^{\frac{s}{2}}(u_n\psi_{\rho})|^2\,{\rm d}x=\lim_{\rho\rightarrow0^{+}}\int_{\mathbb{R}^3}\psi_{\rho}^2{\rm d}\mu=\mu(\{x_j\})
\end{equation*}
and the fist conclusion of $(iii)$ is established.

Note that using H\"{o}lder's inequality, we get that
\begin{align*}
|A_3|&=\Big|2\int_{\mathbb{R}^3}u_n\psi_{\rho}(-\Delta)^{\frac{s}{2}}u_n(-\Delta)^{\frac{s}{2}}\psi_{\rho}\,{\rm d}x\Big|\leq2\Big(\int_{\mathbb{R}^3}u_n^2|(-\Delta)^{\frac{s}{2}}\psi_{\rho}|^2\,{\rm d}x\Big)^{\frac{1}{2}}\Big(\int_{\mathbb{R}^3}\psi_{\rho}^2|(-\Delta)^{\frac{s}{2}}u_n|^2\,{\rm d}x\Big)^{\frac{1}{2}}\\
&\leq C|A_1|^{\frac{1}{2}},
\end{align*}
\begin{align*}
|A_4|&=\Big|-2\int_{\mathbb{R}^3}\psi_{\rho}(-\Delta)^{\frac{s}{2}}u_nB(u_n,\psi_{\rho})\,{\rm d}x\Big|\leq2\Big(\int_{\mathbb{R}^3}B^2(u_n,\psi_{\rho})\,{\rm d}x\Big)^{\frac{1}{2}}\Big(\int_{\mathbb{R}^3}\psi_{\rho}^2|(-\Delta)^{\frac{s}{2}}u_n|^2\,{\rm d}x\Big)^{\frac{1}{2}}\\
&\leq C|A_2|^{\frac{1}{2}}
\end{align*}
and
\begin{align*}
|A_5|&=\Big|-2\int_{\mathbb{R}^3}u_n(-\Delta)^{\frac{s}{2}}\psi_{\rho}B(u_n,\psi_{\rho})\,{\rm d}x\Big|\leq2\Big(\int_{\mathbb{R}^3}u_n^2|(-\Delta)^{\frac{s}{2}}\psi_{\rho}|^2\,{\rm d}x\Big)^{\frac{1}{2}}\Big(\int_{\mathbb{R}^3}B^2(u_n,\psi_{\rho})\,{\rm d}x\Big)^{\frac{1}{2}}\\
&\leq 2|A_1|^{\frac{1}{2}}|A_2|^{\frac{1}{2}}.
\end{align*}
Hence, we only need to show that
\begin{equation}\label{equ2-1}
\lim_{\rho\rightarrow0^{+}}\limsup_{n\rightarrow\infty}\int_{\mathbb{R}^3}u_n^2|(-\Delta)^{\frac{s}{2}}\psi_{\rho}|^2\,{\rm d}x=0
\end{equation}
and
\begin{equation}\label{equ2-2}
\lim_{\rho\rightarrow0^{+}}\limsup_{n\rightarrow\infty}\int_{\mathbb{R}^3}B^2(u_n,\psi_{\rho})(x)\,{\rm d}x=0.
\end{equation}

In fact, using the assumption that $u_n\rightarrow u$ in $L^2(\mathbb{R}^3)$, similar arguments as Lemma 2.8 and 2.9 in \cite{BCSS}, we can conclude that \eqref{equ2-1} and \eqref{equ2-2} hold true.

2. Let $\phi\in C^{\infty}(\mathbb{R}^3)$ such that $\phi=0$ on $B_1(0)$, $\phi=1$ on $\mathbb{R}^3\backslash B_2(0)$, $1\leq\phi\leq1$ and $|\nabla\phi|\leq C$. Set $\phi_R(x)=\phi(\frac{x}{R})$. It follows from Sobolev inequality that
\begin{equation*}
\int_{\mathbb{R}^3}|u_n\phi_R|^{2_s^{\ast}}\,{\rm d}x\leq\Big(\mathcal{S}_s^{-1}\int_{\mathbb{R}^3}|(-\Delta)^{\frac{s}{2}}(u_n\phi_R)|^2\,{\rm d}y\,{\rm d}x\Big)^{\frac{2_s^{\ast}}{2}}.
\end{equation*}
It is easy to check that
\begin{equation*}
\lim_{R\rightarrow+\infty}\limsup_{n\rightarrow\infty}\int_{\mathbb{R}^3}|u_n\phi_R|^{2_s^{\ast}}\,{\rm d}x=\nu_{\infty}.
\end{equation*}
Similar argument to the proof of the first conclusion, we only need to show that
\begin{equation}\label{equ2-3}
\lim_{R\rightarrow+\infty}\limsup_{n\rightarrow\infty}\int_{\mathbb{R}^3}u_n^2|(-\Delta)^{\frac{s}{2}}\phi_R|^2\,{\rm d}x=0
\end{equation}
and
\begin{equation}\label{equ2-4}
\lim_{R\rightarrow+\infty}\limsup_{n\rightarrow\infty}\int_{\mathbb{R}^3}B^2(u_n,\phi_R)\,{\rm d}x=0.
\end{equation}
$\bullet$ Estimate of $(-\Delta)^{\frac{s}{2}}\phi_R$.

\begin{align}\label{equ2-5}
|(-\Delta)^{\frac{s}{2}}\phi_R(x)|&=C_s\Big|{\rm P.V.}\int_{\mathbb{R}^3}\frac{\phi_R(x)-\phi_R(y)}{|x-y|^{3+s}}\,{\rm d}y\Big|\leq C_s{\rm P.V.}\int_{\{|x-y|\leq R\}}\frac{|\nabla \phi_R(\xi)|}{|x-y|^{2+s}}\,{\rm d}y\nonumber\\
&+2C_s\int_{\{|x-y|\geq R\}}\frac{1}{|x-y|^{3+s}}\,{\rm d}y\nonumber\\
&\leq\frac{C}{\rho}{\rm P.V.}\int_{\{|x-y|\leq R\}}\frac{1}{|x-y|^{2+s}}\,{\rm d}y+2C_s\int_{\{|x-y|\geq R\}}\frac{1}{|x-y|^{3+s}}\,{\rm d}y\nonumber\\
&\leq \frac{C}{R^s}.
\end{align}
where $\xi=y+\tau(x-y)$ with $\tau\in(0,1)$.

$\bullet$ Estimate of $I(u_n,\phi_R)$. By H\"{o}lder's inequality, we have that
\begin{align}\label{equ2-6}
&\int_{\mathbb{R}^3}\frac{(\phi_R(x)-\phi_R(y))(u_n(x)-u_n(y))}{|x-y|^{3+s}}\,{\rm d}y\nonumber\\
&\leq\Big(\int_{\mathbb{R}^3}\frac{|\phi_R(x)-\phi_R(y)|^2}{|x-y|^{3+s}}\,{\rm d}y\Big)^{\frac{1}{2}}\Big(\int_{\mathbb{R}^3}\frac{|u_n(x)-u_n(y)|^2}{|x-y|^{3+s}}\,{\rm d}y\Big)^{\frac{1}{2}}\nonumber\\
&\leq2\Big(\int_{\mathbb{R}^3}\frac{|\phi_R(x)-\phi_R(y)|}{|z-y|^{3+s}}\,{\rm d}y\Big)^{\frac{1}{2}}\Big(\int_{\mathbb{R}^3}\frac{|u_n(x)-u_n(y)|^2}{|x-y|^{3+s}}\,{\rm d}y\Big)^{\frac{1}{2}}\nonumber\\
&\leq\frac{C}{R^{\frac{s}{2}}}\Big(\int_{\mathbb{R}^3}\frac{|u_n(x)-u_n(y)|^2}{|x-y|^{3+s}}\,{\rm d}y\Big)^{\frac{1}{2}}.
\end{align}

Therefore, by \eqref{equ2-5}, one has
\begin{align*}
\int_{\mathbb{R}^3}u_n^2|(-\Delta)^{\frac{s}{2}}\phi_R|^2\,{\rm d}x\leq\frac{C}{R^{2s}}\int_{\mathbb{R}^3}u_n^2\,{\rm d}x\leq\frac{C}{R^{2s}}
\end{align*}
which implies that \eqref{equ2-3} holds. From \eqref{equ2-6} and Proposition 3.4 in \cite{NPV}, we deduce that
\begin{align*}
\int_{\mathbb{R}^3}B^2(u_n,\phi_R)\,{\rm d}x&\leq\frac{C}{R^s}\int_{\mathbb{R}^3}\int_{\mathbb{R}^3}\frac{|u_n(x)-u_n(y)|^2}{|x-y|^{3+s}}\,{\rm d}y\,{\rm d}x\\
&=\frac{C}{R^s}\int_{\mathbb{R}^3}|\xi|^s|\widehat{u}_n(\xi)|^2\,{\rm d}\xi\leq\frac{C}{R^s}\int_{\mathbb{R}^3}(1+|\xi|^{2s})|\widehat{u}_n(\xi)|^2\,{\rm d}\xi\\
&=\frac{C}{R^s}\|u_n\|^2\leq\frac{C}{R^s}
\end{align*}
which yields \eqref{equ2-4}. Thus we complete the proof of the second conclusion.
\end{proof}

\section{The penalization problem}

For the bounded domain $\Lambda$ given in $(V_1)$, $k>2$, $a>0$ such that $f(a)+a^{2_s^{\ast}-1}=\frac{\alpha_0}{k}a$ where $\alpha_0$ is mentioned in $(V_0)$, we consider a new problem
\begin{equation}\label{equ3-1}
\left\{
  \begin{array}{ll}
   (-\Delta)^su+V(\varepsilon z)u+\phi u=g(\varepsilon z,u) & \hbox{in $\mathbb{R}^3$,} \\
    (-\Delta)^t\phi=u^2& \hbox{in $\mathbb{R}^3$}
  \end{array}
\right.
\end{equation}
where $g(z,\tau)=\chi( z)(f(\tau)+(\tau^{+})^{2_s^{\ast}-1})+(1-\chi(z))\tilde{f}(\tau)$ with
\begin{equation*}
\tilde{f}(\tau)=\left\{
  \begin{array}{ll}
    f(\tau)+(\tau^{+})^{2_s^{\ast}-1} & \hbox{if $\tau\leq a$,} \\
    \frac{\alpha_0}{k}\tau & \hbox{if $\tau>a$}
  \end{array}
\right.
\end{equation*}
and $\chi(z)$ is a smooth function such that $\chi(z)=1$ on $\Lambda$, $0\leq\chi(z)\leq1$ on $\Lambda\backslash\Lambda'$, $\chi(z)=0$ on $\mathbb{R}^3\backslash\Lambda'$, where $\Lambda'$ is a suitable open set satisfying $\bar{\Lambda}\subset\Lambda'$ and $V(z)>\inf_{\xi\in\Lambda}V(\xi)$ for all $z\in\bar{\Lambda'}\backslash\Lambda$. It is easy to see that under the assumptions $(f_1)$-$(f_3)$, $g(z,\tau)$ is a Caratheodory function and satisfies the following assumptions:\\
$(g_1)$ $g(z,\tau)=o(\tau^3)$ as $\tau\rightarrow0$ uniformly on $z\in\mathbb{R}^3$;\\
$(g_2)$ $g(z,\tau)\leq f(\tau)+\tau^{2_s^{\ast}-1}$ for all $\tau\in\mathbb{R}^{+}$ and $z\in\mathbb{R}^3$, $g(z,\tau)=0$ for all $z\in\mathbb{R}^3$ and $\tau<0$;\\
$(g_3)$ $0<2\tilde{F}(\tau)\leq\tilde{f}(\tau)\tau\leq\frac{\alpha_0}{k}\tau^2\leq\frac{V(x)}{k}\tau^2$ for all $s\geq0$ with the number $k>2$, where $\tilde{F}(\tau)$ is a prime function of $\tilde{f}$;\\
$(g_4)$ $0<4G(z,\tau)\leq g(z,\tau)\tau$ for all $z\in\Lambda$, $\tau>0$ or $z\in\mathbb{R}^3\backslash\Lambda$, $\tau\leq a$ and $g(z,\tau)\tau+\frac{V(z)}{4}\tau^2\geq 4G(z,\tau)>0$ for all $z\in\mathbb{R}^3$ and $\tau>0$, where $G(z,\tau)$ is a prime function of $g(z,\tau)$;\\
$(g_5)$ $\frac{g(z,s\tau)}{s}$ is nondecreasing in $\tau\in\mathbb{R}^{+}$ uniformly for $z\in\mathbb{R}^3$, $\frac{g(z,s\tau)}{\tau^3}$ is nondecreasing in $\tau\in\mathbb{R}^{+}$ and $z\in\Lambda$, $\frac{g(z,s\tau)}{\tau^3}$ is nondecreasing in $\tau\in(0,a)$ and $z\in\mathbb{R}^3\backslash\Lambda'$.

The energy functional corresponding to \eqref{equ3-1} is defined as
\begin{equation*}
J_{\varepsilon}(u)=\frac{1}{2}\int_{\mathbb{R}^3}\Big(|(-\Delta)^{\frac{s}{2}}u|^2+V(\varepsilon z)u^2\Big)\,{\rm d}z+\frac{1}{4}\int_{\mathbb{R}^3}\phi_u^tu^2\,{\rm d}z-\int_{\mathbb{R}^3}G(\varepsilon z,u)\,{\rm d}z, \,\, u\in H_{\varepsilon}
\end{equation*}
and $J_{\varepsilon}\in C^1(H_{\varepsilon},\mathbb{R})$.

By standard argument, the functional $J_{\varepsilon}$ satisfies the mountain pass geometry.
\begin{lemma}\label{lem3-1}
Suppose $(V_0)$, $(V_1)$ and $(f_0)-(f_2)$ hold, then the functional $J_{\varepsilon}$ has the following properties:\\
$(i)$ there exist $\alpha,\rho>0$ such that $J_{\varepsilon}(u)\geq\alpha$ for $\|u\|_{H_\varepsilon}=\rho$;\\
$(ii)$ there exists $e_0\in H_{\varepsilon}$ satisfying $\|e_0\|_{H_\varepsilon}>\rho$ such that $J_{\varepsilon}(e_0)<0$.
\end{lemma}
By Lemma \ref{lem3-1} and Theorem 1.15 in \cite{Willem} (Mountain pass theorem without Palais-Smale condition), it follows that there exists a $(PS)_{c_{\varepsilon}}$ sequence $\{u_n\}\subset H_{\varepsilon}$ such that
\begin{equation}\label{equ3-2}
J_{\varepsilon}(u_n)\rightarrow c_{\varepsilon}\quad\text{and}\quad J'_{\varepsilon}(u_n)\rightarrow0\,\,\text{as}\,\, n\rightarrow\infty,
\end{equation}
where $c_{\varepsilon}=\inf\limits_{\gamma\in \Gamma}\max\limits_{t\in[0,1]}J_{\varepsilon}(\gamma(t))>0$. Here
\begin{equation*}
\Gamma=\{\gamma\in C([0,1],H_{\varepsilon})\,\,|\,\, \gamma(0)=0, J_{\varepsilon}(\gamma(1))<0\}.
\end{equation*}
Similarly to the argument in \cite{PF,Willem},  by $(g_5)$, the equivalent characterization of $c_{\varepsilon}$ is given by
\begin{equation*}
c_{\varepsilon}=\inf_{u\in H_{\varepsilon}\backslash\{0\}}\max_{t\geq0}J_{\varepsilon}(tu)=\inf_{u\in \mathcal{N}_{\varepsilon}}J_{\varepsilon}(u),
\end{equation*}
where $\mathcal{N}_{\varepsilon}$ is the Nehari manifold defined as
\begin{equation*}
\mathcal{N}_{\varepsilon}=\{u\in H_{\varepsilon}\backslash\{0\}\,\,\Big|\,\, \langle J_{\varepsilon}'(u),u\rangle=0\}.
\end{equation*}
For author's convenience, we give the rough proof. We state it as the following Proposition.
\begin{proposition}
\begin{equation*}
c_{\varepsilon}=\inf_{u\in H_{\varepsilon}\backslash\{0\}}\max_{t\geq0}J_{\varepsilon}(tu)=\inf_{u\in \mathcal{N}_{\varepsilon}}J_{\varepsilon}(u).
\end{equation*}
\end{proposition}
\begin{proof}
1. For each $u\in H_{\varepsilon}\backslash\{0\}$, we claim that there exists a unique $t_u>0$ such that $t_uu\in\mathcal{N}_{\varepsilon}$. Indeed, set $h(t)=J_{\varepsilon}(tu)$, by $(g_1)$ and $(g_2)$, it is easy to check that $h(t)>0$ when $t>0$ small and $h(t)<0$ when $t>0$ large. Since $h\in C^1(\mathbb{R}^{+},\mathbb{R})$ and $h(0)=0$, there is $t_u>0$ global maximum point of $h(t)$ and $h'(t)=0$. Thus, $\langle J_{\varepsilon}(t_u u),t_uu\rangle=0$, and $t_uu\in\mathcal{N}_{\varepsilon}$. We see that $t_u>0$ is the unique positive number such that $h'(t_u)=0$. Indeed, suppose by contradiction that there exist $t_1>t_2>0$ such that $h'(t_1)=h'(t_2)=0$. Then for $i=1,2$,
\begin{equation*}
t_i\|u\|_{H_{\varepsilon}}^2+t_i^3\int_{\mathbb{R}^3}\phi_u^tu^2\,{\rm d}z=\int_{\mathbb{R}^3}g(\varepsilon z,t_iu)u\,{\rm d}z.
\end{equation*}
Therefore
\begin{align*}
(\frac{1}{t_1^2}-\frac{1}{t_2^2})\|u\|_{H_{\varepsilon}}^2=\int_{\mathbb{R}^3}\Big(\frac{g(\varepsilon z, t_1u)}{(t_1u)^3}-\frac{g(\varepsilon z, t_2u)}{(t_2u)^3}\Big)u^4\,{\rm d}z
\end{align*}

$\bullet$ If ${\rm supp}(u^{+})\subset\Lambda/\varepsilon$, then $g(\varepsilon z,u)=f(u)+(u^{+})^{2_s^{\ast}-1}$, the uniqueness of $t_u$ follows from the hypothesis $(f_2)$.

$\bullet$ If ${\rm supp}(u^{+})\subset\mathbb{R}^3\backslash(\Lambda'/\varepsilon)$, then $g(\varepsilon z,u)=\tilde{f}(u)$. By the definition of $\tilde{f}$, we have that
\begin{align*}
(\frac{1}{t_1^2}-\frac{1}{t_2^2})\|u\|_{H_{\varepsilon}}^2&\geq\int_{\mathbb{R}^3\backslash(\Lambda'/\varepsilon)\cap\{t_2u<a<t_1u\}}\Big(\frac{g(\varepsilon z, t_1u)}{(t_1u)^3}-\frac{g(\varepsilon z, t_2u)}{(t_2u)^3}\Big)u^4\,{\rm d}z\\
&+\int_{\mathbb{R}^3\backslash(\Lambda'/\varepsilon)\cap\{a<t_2u\}}\Big(\frac{g(\varepsilon z, t_1u)}{(t_1u)^3}-\frac{g(\varepsilon z, t_2u)}{(t_2u)^3}\Big)u^4\,{\rm d}z\\
&=\int_{\mathbb{R}^3\backslash(\Lambda'/\varepsilon)\cap\{t_2u<a<t_1u\}}\Big(\frac{\alpha_0}{k}\frac{u^2}{t_1^2}-\frac{f(t_2u)+t_2^{2_s^{\ast}-1}(u^{+})^{2_s^{\ast}-1}}{(t_2u)^3}u^4\Big)\,{\rm d}z\\
&+\int_{\mathbb{R}^3\backslash(\Lambda'/\varepsilon)\cap\{a<t_2u\}}\frac{\alpha_0}{k}\Big(\frac{1}{t_1^2}-\frac{1}{t_2^2}\Big)u^2\,{\rm d}z.
\end{align*}
Multiplying both sides by $\frac{1}{\frac{1}{t_1^2}-\frac{1}{t_2^2}}$ and using the hypothesis $t_1>t_2>0$, we get
\begin{align*}
\|u\|_{H_{\varepsilon}}^2&\leq\frac{1}{k}\int_{\mathbb{R}^3\backslash(\Lambda'/\varepsilon)\cap\{a<t_2u\}}\alpha_0u^2\,{\rm d}z-\frac{t_2^2}{t_1^2-t_2^2}\int_{\mathbb{R}^3\backslash(\Lambda'/\varepsilon)\cap\{t_2u<a<t_1u\}}\frac{\alpha_0}{k}u^2\,{\rm d}z\\
&+\frac{t_1^2}{t_1^2-t_2^2}\int_{\mathbb{R}^3\backslash(\Lambda'/\varepsilon)\cap\{t_2u<a<t_1u\}}\frac{f(t_2u)+(t_2u^{+})^{2_s^{\ast}-1}}{t_2}u\,{\rm d}z\\
&\leq\frac{1}{k}\int_{\mathbb{R}^3\backslash(\Lambda'/\varepsilon)\cap\{a<t_2u\}}\alpha_0u^2\,{\rm d}z-\frac{t_2^2}{t_1^2-t_2^2}\int_{\mathbb{R}^3\backslash(\Lambda'/\varepsilon)\cap\{t_2u<a<t_1u\}}\frac{\alpha_0}{k}u^2\,{\rm d}z\\
&+\frac{t_1^2}{t_1^2-t_2^2}\int_{\mathbb{R}^3\backslash(\Lambda'/\varepsilon)\cap\{t_2u<a<t_1u\}}\frac{\alpha_0}{k}u^2\,{\rm d}z\\
&\leq\frac{1}{k}\int_{\mathbb{R}^3}\alpha_0u^2\,{\rm d}z\leq\frac{1}{k}\|u\|_{H_{\varepsilon}}^2.
\end{align*}
Since $u\neq0$, we have that $k\leq1$, but this is a contradiction. Thus, the uniqueness of $t_u$ follows for ${\rm supp}(u^{+})\subset\mathbb{R}^3\backslash(\Lambda'/\varepsilon)$.

$\bullet$ If ${\rm supp}(u^{+})\subset(\mathbb{R}^3\backslash(\Lambda/\varepsilon))\backslash(\mathbb{R}^3\backslash(\Lambda'/\varepsilon))$, by the definition of $g$ and hypothesis $(f_2)$, we have that
\begin{align*}
&(\frac{1}{t_1^2}-\frac{1}{t_2^2})\|u\|_{H_{\varepsilon}}^2=\int_{(\mathbb{R}^3\backslash(\Lambda/\varepsilon))\backslash(\mathbb{R}^3\backslash(\Lambda'/\varepsilon))}(1-\chi(\varepsilon z))\Big(\frac{\tilde{f}(t_1 u)}{(t_1u)^3}-\frac{\tilde{f}(t_2 u)}{(t_2u)^3}\Big)u^4\,{\rm d}z\\
&+\int_{(\mathbb{R}^3\backslash(\Lambda/\varepsilon))\backslash(\mathbb{R}^3\backslash(\Lambda'/\varepsilon))}\chi(\varepsilon z)\Big(\frac{f(t_1u)+(t_1u^{+})^{2_s^{\ast}-1}}{(t_1u)^3}-\frac{(f(t_2u)+(t_2u^{+})^{2_s^{\ast}-1})}{(t_2u)^3}\Big)u^4\,{\rm d}z\\
&\geq\int_{(\mathbb{R}^3\backslash(\Lambda/\varepsilon))\backslash(\mathbb{R}^3\backslash(\Lambda'/\varepsilon))}(1-\chi(\varepsilon z))\Big(\frac{\tilde{f}(t_1 u)}{(t_1u)^3}-\frac{\tilde{f}(t_2 u)}{(t_2u)^3}\Big)u^4\,{\rm d}z\\
&\geq\int_{(\mathbb{R}^3\backslash(\Lambda/\varepsilon))\backslash(\mathbb{R}^3\backslash(\Lambda'/\varepsilon))\cap\{t_2u<a<t_1u\}}(1-\chi(\varepsilon z))\Big(\frac{1}{t_1^2}\frac{\alpha_0}{k}u^2-\frac{(f(t_2u)+(t_2u^{+})^{2_s^{\ast}-1})}{(t_2u)^3}\Big)u^4\,{\rm d}z\\
&+\Big(\frac{1}{t_1^2}-\frac{1}{t_2^2}\Big)\frac{1}{k}\int_{(\mathbb{R}^3\backslash(\Lambda/\varepsilon))\backslash(\mathbb{R}^3\backslash(\Lambda'/\varepsilon))\cap\{a<t_2u\}}(1-\chi(\varepsilon z))\alpha_0u^2\,{\rm d}z.
\end{align*}
Multiplying both sides by $\frac{1}{\frac{1}{t_1^2}-\frac{1}{t_2^2}}$ and using the hypothesis $t_1>t_2>0$, we get
\begin{align*}
\|u\|_{H_{\varepsilon}}^2&\leq\frac{1}{k}\int_{(\mathbb{R}^3\backslash(\Lambda/\varepsilon))\backslash(\mathbb{R}^3\backslash(\Lambda'/\varepsilon))\cap\{a<t_2u\}}(1-\chi(\varepsilon z))\alpha_0u^2\,{\rm d}z\\
&-\frac{t_2^2}{t_1^2-t_2^2}\int_{(\mathbb{R}^3\backslash(\Lambda/\varepsilon))\backslash(\mathbb{R}^3\backslash(\Lambda'/\varepsilon))\cap\{t_2u<a<t_1u\}}(1-\chi(\varepsilon z))\frac{\alpha_0}{k}u^2\,{\rm d}z\\
&+\frac{t_1^2}{t_1^2-t_2^2}\int_{(\mathbb{R}^3\backslash(\Lambda/\varepsilon))\backslash(\mathbb{R}^3\backslash(\Lambda'/\varepsilon))\cap\{t_2u<a<t_1u\}}(1-\chi(\varepsilon z))\frac{\alpha_0}{k}u^2\,{\rm d}z\\
&\leq\frac{1}{k}\int_{\mathbb{R}^3}(1-\chi(\varepsilon z))\alpha_0u^2\,{\rm d}z\leq\frac{1}{k}\int_{\mathbb{R}^3}\alpha_0u^2\,{\rm d}z\leq\frac{1}{k}\|u\|_{H_{\varepsilon}}^2
\end{align*}
which implies a contradiction with $k>2$.

$\bullet$ If ${\rm supp}(u^{+})\subset(\mathbb{R}^3\backslash(\Lambda/\varepsilon))\backslash(\mathbb{R}^3\backslash(\Lambda'/\varepsilon))\cup\mathbb{R}^3\backslash(\Lambda'/\varepsilon)$, similar argument as the above, we can deduce that $k\leq2$, this is a contradiction.

Therefore, whatever any cases, the claim holds true. Thus
\begin{equation*}
\inf_{u\in H_{\varepsilon}\backslash\{0\}}\max_{t\geq0}J_{\varepsilon}(tu)=\inf_{u\in \mathcal{N}_{\varepsilon}}J_{\varepsilon}(u).
\end{equation*}

By $(ii)$ of Lemma \ref{lem3-1}, using standard argument, we get that
\begin{equation*}
\inf_{u\in H_{\varepsilon}\backslash\{0\}}\max_{t\geq0}J_{\varepsilon}(tu)\geq\inf_{\gamma\in\Gamma}\max_{t\in[0,1]}J_{\varepsilon}(\gamma(t)).
\end{equation*}

For any $\gamma\in\Gamma$, $\gamma([0,1])\cap\mathcal{N}_{\varepsilon}\neq\emptyset$. Indeed, if $u\in H_{\varepsilon}\backslash\{0\}$ is interior to or on $\mathcal{N}_{\varepsilon}$, then
\begin{align*}
\|u\|_{H_{\varepsilon}}^2+\int_{\mathbb{R}^3}\phi_u^tu^2\,{\rm d}z\geq\int_{\mathbb{R}^3}g(\varepsilon z,u)u\,{\rm d}z
\end{align*}
and
\begin{align*}
4J_{\varepsilon}(u)=\langle J_{\varepsilon}'(u),u\rangle+\|u\|_{H_{\varepsilon}}+\int_{\mathbb{R}^3}\Big(g(\varepsilon z,u)-4G(\varepsilon z,u)\Big)\,{\rm d}z>0.
\end{align*}
Hence $\gamma$ crosses $\mathcal{N}_{\varepsilon}$ since $\gamma(0)=0$, $J_{\varepsilon}(\gamma(1))\leq0$ and $\gamma(1)\neq0$. Therefore
\begin{equation*}
\max_{t\in[0,1]}J_{\varepsilon}(\gamma(t))\geq\inf_{\mathcal{N}_{\varepsilon}}J_{\varepsilon}(u).
\end{equation*}
\end{proof}

The following Lemma gives the estimate of the critical value $c_{\varepsilon}$.
\begin{lemma}\label{lem3-2}
Suppose that $(V_0)$, $(V_1)$ and $(f_0)-(f_2)$ hold, then the infinimum $c_{\varepsilon}$ satisfies
\begin{equation*}
0<c_{\varepsilon}<\frac{s}{3}\mathcal{S}_s^{\frac{3}{2s}}
\end{equation*}
for $\varepsilon$ small enough, where $\mathcal{S}_s$ is the best Sobolev constant for the embedding $\mathcal{D}^{s,2}(\mathbb{R}^3)\hookrightarrow L^{2_s^{\ast}}(\mathbb{R}^3)$.
\end{lemma}

\begin{proof}
Without loss of generalization, we assume that $0\in\Lambda$. Choose $R>0$ such that $B_{2R}(0)\subset \Lambda/\varepsilon$ and $\psi\in C_0^{\infty}(B_{2R}(0))$ satisfying $\psi=1$ on $B_R(0)$ and $0\leq\psi\leq1$ on $B_{2R}(0)$. Given $\varepsilon>0$, we define
\begin{equation*}
v_{\varepsilon}(z)=\psi(z)U_{\varepsilon}(z),\quad x\in\mathbb{R}^3,
\end{equation*}
where $U_{\varepsilon}(z)=\varepsilon^{-\frac{3-2s}{2}}u^{\ast}(z/\varepsilon)$, $u^{\ast}(z)=\frac{\widetilde{u}(z/\mathcal{S}_s^{\frac{1}{2s}})}{\|\widetilde{u}\|_{2_s^{\ast}}}$, $\kappa\in\mathbb{R}\backslash\{0\}$, $\mu>0$ and $x_0\in\mathbb{R}^3$ are fixed constants, $\widetilde{u}(z)=\kappa(\mu^2+|z-x_0|^2)^{-\frac{3-2s}{2}}$. From Proposition 21 and Proposition 22 in \cite{SV}, Lemma 3.3 in \cite{Teng}, we know that
\begin{equation}\label{equ3-3}
\int_{\mathbb{R}^{3}}|(-\Delta)^{\frac{s}{2}}v_{\varepsilon}|^2\,{\rm d}z\leq\mathcal{S}_s^{\frac{3}{2s}}+O(\varepsilon^{3-2s}),
\end{equation}
\begin{equation}\label{equ3-4}
\int_{\mathbb{R}^{3}}|v_{\varepsilon}|^{2_s^{\ast}}\,{\rm d}z=\mathcal{S}_s^{\frac{3}{2s}}+O(\varepsilon^3),
\end{equation}
and
\begin{equation}\label{equ3-5}
\int_{\mathbb{R}^3}|v_{\varepsilon}(z)|^p\,{\rm d}z=\left\{
\begin{array}{ll}
O(\varepsilon^{\frac{(2-p)3+2sp}{2}}),&\hbox{$p>\frac{3}{3-2s}$,} \\
O(\varepsilon^{\frac{(2-p)3+2sp}{2}}|\log\varepsilon|), & \hbox{$p=\frac{3}{3-2s}$,} \\
O(\varepsilon^{\frac{3-2s}{2}p}), & \hbox{$1< p<\frac{3}{3-2s}$.}
\end{array}
\right.
\end{equation}

Since ${\rm supp}(v_{\varepsilon})\subset \Lambda/\varepsilon$, $g(\varepsilon z,v_{\varepsilon})=f(v_{\varepsilon})+v_{\varepsilon}^{2_s^{\ast}-1}$. There exists $t_{\varepsilon}>0$ such that $\sup\limits_{t\geq0}J_{\varepsilon}(tv_{\varepsilon})=J_{\varepsilon}(t_{\varepsilon}v_{\varepsilon})$. Hence $\frac{{\rm d}J_{\varepsilon}(tv_{\varepsilon})}{{\rm d} t}\Big|_{t=t_{\varepsilon}}=0$, that is
\begin{align*}
t_{\varepsilon}\int_{\mathbb{R}^3}(|(-\Delta)^{\frac{s}{2}}v_{\varepsilon}|^2+V(\varepsilon z) v_{\varepsilon}^2)\,{\rm d}z+t_{\varepsilon}^3\int_{\mathbb{R}^3}\phi_{v_{\varepsilon}}^tv_{\varepsilon}^2\,{\rm d}z=&\int_{\mathbb{R}^3} f(t_{\varepsilon}v_{\varepsilon})v_{\varepsilon}\,{\rm d}z\\
&+t_{\varepsilon}^{2_s^{\ast}-1}\int_{\mathbb{R}^3}v_{\varepsilon}^{2_s^{\ast}}\,{\rm d}z.
\end{align*}
By $(f_0)$, we have that
\begin{align*}
\int_{\mathbb{R}^3}(|(-\Delta)^{\frac{s}{2}}v_{\varepsilon}|^2+V(\varepsilon z) v_{\varepsilon}^2)\,{\rm d}z+t_{\varepsilon}^2\int_{\mathbb{R}^3}\phi_{v_{\varepsilon}}^tv_{\varepsilon}^2\,{\rm d}z\geq t_{\varepsilon}^{2_s^{\ast}-2}\int_{\mathbb{R}^3}v_{\varepsilon}^{2_s^{\ast}}\,{\rm d}z.
\end{align*}
It follows from $(iii)$ of Lemma \ref{lem2-1} that
\begin{align}\label{equ3-6}
t_{\varepsilon}^{2_s^{\ast}-2}\leq\frac{1}{\|v_{\varepsilon}\|_{2_s^{\ast}}^{2_s^{\ast}}}\Big(\int_{\mathbb{R}^3}(|(-\Delta)^{\frac{s}{2}}v_{\varepsilon}|^2+V(\varepsilon z) v_{\varepsilon}^2)\,{\rm d}z+Ct_{\varepsilon}^2\|v_{\varepsilon}\|_{\frac{12}{3+2t}}^4\Big).
\end{align}
Thus, \eqref{equ3-3}-\eqref{equ3-6} imply that $t_{\varepsilon}\leq C_1$, where $C_1$ is independent of $\varepsilon>0$ small. On the other hand, we may assume that there is a positive constant $C_2>0$ such that $t_{\varepsilon}\geq C_2>0$ for $\varepsilon>0$ small. Otherwise, we can find a sequence $\varepsilon_n\rightarrow0$ as $n\rightarrow\infty$ such that $t_{\varepsilon_n}\rightarrow0$ as $n\rightarrow\infty$. Therefore
\begin{equation*}
0<c_{\varepsilon}\leq \sup_{t\geq0}J_{\varepsilon}(tv_{\varepsilon_n})=J_{\varepsilon}(t_{\varepsilon_n}v_{\varepsilon_n})\rightarrow0,
\end{equation*}
which is a contradiction.

Denote $g(t)=\frac{t^2}{2}\int_{\mathbb{R}^3}|(-\Delta)^{\frac{s}{2}}v_{\varepsilon}|^2\,{\rm d}z-\frac{ t^{2_s^{\ast}}}{2_s^{\ast}}\int_{\mathbb{R}^3}|v_{\varepsilon}|^{2_s^{\ast}}\,{\rm d}z$, by \eqref{equ3-3} and \eqref{equ3-4}, it is easy to check that
\begin{align*}
\sup_{t\geq0}g(t)=\frac{s}{3}\frac{\Big(\int_{\mathbb{R}^3}|(-\Delta)^{\frac{s}{2}}v_{\varepsilon}|^2\,{\rm d}z\Big)^{\frac{3}{2s}}}{\Big(\int_{\mathbb{R}^3}|v_{\varepsilon}|^{2_s^{\ast}}\,{\rm d}z\Big)^{\frac{3-2s}{2s}}}&\leq\frac{s}{3}\frac{\Big(\mathcal{S}_s^{\frac{3}{2s}}+O(\varepsilon^{3-2s})\Big)^{\frac{3}{2s}}}{\Big(\mathcal{S}_s^{\frac{3}{2s}}+O(\varepsilon^3)\Big)^{\frac{3-2s}{2s}}}\\
&=\frac{s}{3}\mathcal{S}_s^{\frac{3}{2s}}+O(\varepsilon^{3-2s}),
\end{align*}
where we have used the elementary inequality $(a+b)^{\alpha}\leq a^{\alpha}+\alpha(a+b)^{\alpha-1}b$, $\alpha\geq1$, $a,b>0$.
Thus, using the fact that $\max\limits_{z\in B_{2R}(0)}V(\varepsilon z)\leq C$ for some $C>0$ independent of $\varepsilon$ and $s>\frac{3}{4}\Rightarrow2<\frac{3}{3-2s}$, by \eqref{equ3-5}, we deduce that
\begin{align}\label{equ3-7}
J_{\varepsilon}(t_{\varepsilon}v_{\varepsilon})&\leq \sup_{t\geq0}g(t)+C\int_{\mathbb{R}^3}V(\varepsilon z)|v_{\varepsilon}|^2\, {\rm d}z+C\int_{\mathbb{R}^3}\phi_{v_{\varepsilon}}^tv_{\varepsilon}^2\,{\rm d}z-C\int_{\mathbb{R}^3}|v_{\varepsilon}|^{q}\, {\rm d}z\nonumber\\
&\leq\frac{s}{3}\mathcal{S}_s^{\frac{3}{2s}}+O(\varepsilon^{3-2s})+C\int_{\mathbb{R}^3}|v_{\varepsilon}|^2\, {\rm d}z+C\Big(\int_{\mathbb{R}^3}|v_{\varepsilon}|^{\frac{12}{3+2t}}\,{\rm d}z\Big)^{\frac{3+2t}{3}}\nonumber\\
&-C\int_{\mathbb{R}^3}|v_{\varepsilon}|^q\, {\rm d}z\nonumber\\
&\leq\frac{s}{3}\mathcal{S}_s^{\frac{3}{2s}}+O(\varepsilon^{3-2s})+C\Big(\int_{\mathbb{R}^3}|v_{\varepsilon}|^{\frac{12}{3+2t}}\,{\rm d}z\Big)^{\frac{3+2t}{3}}-C\int_{\mathbb{R}^3}|v_{\varepsilon}|^q\, {\rm d}z.
\end{align}

By \eqref{equ3-5}, we have that
\begin{align}\label{equ3-8}
\lim_{\varepsilon\rightarrow0^{+}}\frac{\Big(\int_{\mathbb{R}^3}|v_{\varepsilon}|^{\frac{12}{3+2t}}\,{\rm d}z\Big)^{\frac{3+2t}{3}}}{\varepsilon^{3-2s}}\leq\left\{
                                                                                                                            \begin{array}{ll}
                                                                                                                              \lim\limits_{\varepsilon\rightarrow0^{+}}\frac{O(\varepsilon^{2t+4s-3})}{\varepsilon^{3-2s}}=0, & \hbox{$\frac{12}{3+2t}>\frac{3}{3-2s}$,} \\
                                                                                                                              \lim\limits_{\varepsilon\rightarrow0^{+}}\frac{O(\varepsilon^{2t+4s-3}|\log\varepsilon|)}{\varepsilon^{3-2s}}=0, & \hbox{$\frac{12}{3+2t}=\frac{3}{3-2s}$,} \\
                                                                                                                             \lim\limits_{\varepsilon\rightarrow0^{+}} \frac{O(\varepsilon^{2(3-2s)})}{\varepsilon^{3-2s}}=0, & \hbox{$\frac{12}{3+2t}<\frac{3}{3-2s}$.}
                                                                                                                            \end{array}
                                                                                                                          \right.
\end{align}
Since $s>\frac{3}{4}$ and $q\geq4$, then $q>\frac{3}{3-2s}$, $2s-\frac{3-2s}{2}q<0$. Thus
\begin{align}\label{equ3-9}
\lim_{\varepsilon\rightarrow0^{+}}\frac{\int_{\mathbb{R}^3}|v_{\varepsilon}|^{q}\,{\rm d}x}{\varepsilon^{3-2s}}=\lim_{\varepsilon\rightarrow0^{+}}\frac{O(\varepsilon^{3-\frac{3-2s}{2}q})}{\varepsilon^{3-2s}}=+\infty.
\end{align}

Therefore, combining with \eqref{equ3-7}, \eqref{equ3-8} and \eqref{equ3-9}, we conclude that
\begin{equation*}
J_{\varepsilon}(t_{\varepsilon}v_{\varepsilon})<\frac{s}{3}\mathcal{S}_s^{\frac{3}{2s}}
\end{equation*}
for $\varepsilon$ small enough and thus the proof is completed.
\end{proof}

Now we study the $(PS)_{c_{\varepsilon}}$ sequence given in \eqref{equ3-2}.
\begin{lemma}\label{lem3-3}
Sequence $\{u_n\}$ given in \eqref{equ3-2} is bounded in $H_{\varepsilon}$.
\end{lemma}
\begin{proof}
By \eqref{equ3-2} and $(g_2)$-$(g_3)$, we have that
\begin{align*}
c_{\varepsilon}+o(1)&=J_{\varepsilon}(u_n)-\frac{1}{4}\langle J_{\varepsilon}'(u_n),u_n\rangle\\
&=\frac{1}{4}\|u_n\|_{H_{\varepsilon}}^2+\frac{1}{4}\int_{\mathbb{R}^3}\Big(g(\varepsilon z,u_n)u_n-4G(\varepsilon z,u_n)\Big)\,{\rm d}z\\
&=\frac{1}{4}\|u_n\|_{H_{\varepsilon}}^2+\frac{1}{4}\Big(\int_{\Lambda/\varepsilon}+\int_{\Lambda'/\varepsilon\backslash\Lambda/\varepsilon}+\int_{\mathbb{R}^3\backslash\Lambda'/\varepsilon}\Big)\Big(g(\varepsilon z,u_n)u_n-4G(\varepsilon z,u_n)\Big)\,{\rm d}z\\
&\geq\frac{1}{4}\|u_n\|_{H_{\varepsilon}}^2+\frac{1}{4}\int_{\Lambda'/\varepsilon\backslash\Lambda/\varepsilon}\Big((1-\chi(\varepsilon z))(\tilde{f}(u_n)u_n-4\tilde{F}(u_n))\\
&+\chi(\varepsilon z)(f(u_n)u_n-4F(u_n))\Big)\,{\rm d}z+\frac{1}{4}\int_{\mathbb{R}^3\backslash\Lambda'/\varepsilon}(\tilde{f}(u_n)u_n-4\tilde{F}(u_n))\,{\rm d}z\\
&\geq\frac{1}{4}\|u_n\|_{H_{\varepsilon}}^2+\frac{1}{4}\int_{\Lambda'/\varepsilon\backslash\Lambda/\varepsilon}(1-\chi(\varepsilon z))(\tilde{f}(u_n)u_n-4\tilde{F}(u_n))\,{\rm d}z\\
&+\frac{1}{4}\int_{\mathbb{R}^3\backslash\Lambda'/\varepsilon}(\tilde{f}(u_n)u_n-4\tilde{F}(u_n))\,{\rm d}z\\
&\geq\frac{1}{4}\|u_n\|_{H_{\varepsilon}}^2-\frac{1}{2}\int_{\Lambda'/\varepsilon\backslash\Lambda/\varepsilon}(1-\chi(\varepsilon z))
\tilde{F}(u_n)\,{\rm d}z\\
&-\frac{1}{2}\int_{\mathbb{R}^3\backslash\Lambda'/\varepsilon}\tilde{F}(u_n)\,{\rm d}z\\
&\geq\frac{1}{4}\|u_n\|_{H_{\varepsilon}}^2-\frac{1}{4k}\int_{\Lambda'/\varepsilon\backslash\Lambda/\varepsilon}(1-\chi(\varepsilon z))
V(\varepsilon z)|u_n|^2\,{\rm d}z\\
&-\frac{1}{4k}\int_{\mathbb{R}^3\backslash\Lambda'/\varepsilon}V(\varepsilon z)|u_n|^2\,{\rm d}z\\
&\geq\frac{1}{4}(1-\frac{1}{k})\|u_n\|_{H_{\varepsilon}}^2.
\end{align*}
By the choice of $k$, we get the boundedness of $\{u_n\}$ in $H_{\varepsilon}$.
\end{proof}
Next, we show the bounded sequence $\{u_n\}$ is nonvanishing, that is
\begin{lemma}\label{lem3-4}
There exist a sequence $\{z_n\}\subset\mathbb{R}^3$ and $R>0$, $\beta>0$ such that
\begin{equation*}
\liminf_{n\rightarrow\infty}\int_{B_R(z_n)}|u_n|^2\,{\rm d}z\geq \beta,
\end{equation*}
where $\{u_n\}$ is the sequence given by Lemma \ref{lem3-3}.
\end{lemma}

\begin{proof}
Suppose by contradiction that the Lemma does not hold. Thus by the vanishing Lemma, it follows that
\begin{equation}\label{equ3-10}
\int_{\mathbb{R}^3}|u_n|^r\,{\rm d}z\rightarrow0\,\,\text{as}\,\, n\rightarrow\infty\,\, \text{for all}\,\, 2<r<2_s^{\ast}.
\end{equation}
From $(f_1)$-$(f_2)$, $(iii)$ of Lemma \ref{lem2-1} and \eqref{equ3-10}, it is easy to check that
\begin{equation*}
\int_{\mathbb{R}^3}F(u_n)\,{\rm d}z\rightarrow0,\quad \int_{\mathbb{R}^3}f(u_n)u_n\,{\rm d}z\rightarrow0\quad \text{as}\,\, n\rightarrow\infty,
\end{equation*}
and
\begin{equation*}
\int_{\mathbb{R}^3}\phi_{u_n}^tu_n^2\,{\rm d}z\rightarrow0\quad\text{as}\,\, n\rightarrow\infty.
\end{equation*}
By the definition $\tilde{f}$ and \eqref{equ3-10}, we have that
\begin{align*}
\int_{\mathbb{R}^3}G(\varepsilon z,u_n)\,{\rm d}z&=\frac{1}{2_s^{\ast}}\int_{\mathbb{R}^3}\Big(\chi(\varepsilon z)(u_n^{+})^{2_s^{\ast}}+(1-\chi(\varepsilon z))\tilde{F}(u_n)\Big)\,{\rm d}z+o(1)\\
&=\frac{1}{2_s^{\ast}}\int_{\mathbb{R}^3}\chi(\varepsilon z)(u_n^{+})^{2_s^{\ast}}\,{\rm d}z+\frac{1}{2_s^{\ast}}\int_{\{u_n<a\}}(1-\chi(\varepsilon z))(u_n^{+})^{2_s^{\ast}}\,{\rm d}z\\
&+\frac{\alpha_0}{2k}\int_{\{u_n>a\}}(1-\chi(\varepsilon z))|u_n|^2\,{\rm d}z+o(1)
\end{align*}
and
\begin{align*}
\int_{\mathbb{R}^3}g(\varepsilon z,u_n)u_n\,{\rm d}z&=\int_{\mathbb{R}^3}\Big(\chi(\varepsilon z)(u_n^{+})^{2_s^{\ast}}+(1-\chi(\varepsilon z))\tilde{f}(u_n)u_n\Big)\,{\rm d}z+o(1)\\
&=\int_{\mathbb{R}^3}\chi(\varepsilon z)(u_n^{+})^{2_s^{\ast}}\,{\rm d}z+\int_{\{u_n<a\}}(1-\chi(\varepsilon z))(u_n^{+})^{2_s^{\ast}}\,{\rm d}z\\
&+\frac{\alpha_0}{k}\int_{\{u_n>a\}}(1-\chi(\varepsilon z))|u_n|^2\,{\rm d}z+o(1).
\end{align*}
From \eqref{equ3-2}, it follows that
\begin{align}\label{equ3-11}
\frac{1}{2}\|u_n\|_{H_{\varepsilon}}^2-\frac{1}{2_s^{\ast}}\int_{\mathbb{R}^3}\chi(\varepsilon z)(u_n^{+})^{2_s^{\ast}}\,{\rm d}z&-\frac{1}{2_s^{\ast}}\int_{\{u_n<a\}}(1-\chi(\varepsilon z))(u_n^{+})^{2_s^{\ast}}\,{\rm d}z\nonumber\\
&-\frac{\alpha_0}{2k}\int_{\{u_n>a\}}(1-\chi(\varepsilon z))|u_n|^2\,{\rm d}z=c_{\varepsilon}+o(1)
\end{align}
and
\begin{align*}
\|u_n\|_{H_{\varepsilon}}^2-\int_{\mathbb{R}^3}\chi(\varepsilon z)(u_n^{+})^{2_s^{\ast}}\,{\rm d}z&-\int_{\{u_n<a\}}(1-\chi(\varepsilon z))(u_n^{+})^{2_s^{\ast}}\,{\rm d}z\\
&-\frac{\alpha_0}{k}\int_{\{u_n>a\}}(1-\chi(\varepsilon z))|u_n|^2\,{\rm d}z=o(1).
\end{align*}
We may assume that
\begin{equation*}
\|u_n\|_{H_{\varepsilon}}^2-\frac{\alpha_0}{k}\int_{\{u_n>a\}}(1-\chi(\varepsilon z))|u_n|^2\,{\rm d}z\rightarrow l
\end{equation*}
and
\begin{equation*}
\int_{\mathbb{R}^3}\chi(\varepsilon z)(u_n^{+})^{2_s^{\ast}}\,{\rm d}z+\int_{\{u_n<a\}}(1-\chi(\varepsilon z))(u_n^{+})^{2_s^{\ast}}\,{\rm d}z\rightarrow l\quad \text{as}\,\, n\rightarrow\infty.
\end{equation*}
Observe that
\begin{equation*}
V(\varepsilon z)-\frac{\alpha_0}{k}(1-\chi(\varepsilon z))\geq V(\varepsilon z)-\frac{\alpha_0}{k}\geq\frac{V(\varepsilon z)}{2}\quad \text{for any} \,\,z\in\mathbb{R}^3,
\end{equation*}
thus it is easy to check that $l>0$, if not, $\|u_n\|_{H_{\varepsilon}}\rightarrow0$ as $n\rightarrow\infty$, which contradicts with $c_{\varepsilon}>0$.

By \eqref{equ3-11}, we get
\begin{equation}\label{equ3-12}
c_{\varepsilon}=\frac{s}{3}l.
\end{equation}
In view of the definition of $\mathcal{S}_s$, we see that
\begin{align*}
\|u_n\|_{H_{\varepsilon}}^2&-\frac{\alpha_0}{k}\int_{\{u_n>a\}}(1-\chi(\varepsilon z))|u_n|^2\,{\rm d}z\geq\mathcal{S}_{s}\Big(\int_{\mathbb{R}^3}(u_n^{+})^{2_s^{\ast}}\,{\rm d}z\Big)^{\frac{3-2s}{3}}\\
&\geq\mathcal{S}_{s}\Big(\int_{\mathbb{R}^3}\chi(\varepsilon z)(u_n^{+})^{2_s^{\ast}}\,{\rm d}z+\int_{\{u_n<a\}}(1-\chi(\varepsilon z))(u_n^{+})^{2_s^{\ast}}\,{\rm d}z\Big)^{\frac{3-2s}{3}}
\end{align*}
which achieves that
\begin{equation*}
c_{\varepsilon}\geq\frac{s}{3}\mathcal{S}_s^{\frac{3}{2s}},
\end{equation*}
which contradicts with Lemma \ref{lem3-2}.
\end{proof}

\begin{lemma}\label{lem3-5}
The sequence $\{z_n\}$ obtained in Lemma \ref{lem3-3} is bounded in $\mathbb{R}^3$.
\end{lemma}

\begin{proof}
For each $\rho>0$, consider a smooth cut-off function $0\leq \psi_{\rho}\leq1$ such that $\psi_{\rho}=1$ on $|z|\geq 2\rho$, $\psi_{\rho}=0$ on $|z|\leq\rho$ and $|\nabla\psi_{\rho}|\leq\frac{C}{\rho}$. Clearly, $\psi_{\rho}u_n\in H_{\varepsilon}$ for each $\rho>0$. Using $\langle J_{\varepsilon}'(u_n), \psi_{\rho}u_n\rangle=o(1)$, we obtain
\begin{align*}
\int_{\mathbb{R}^3}(-\Delta)^{\frac{s}{2}}u_n(-\Delta)^{\frac{s}{2}}(\psi_{\rho}u_n)\,{\rm d}z+\int_{\mathbb{R}^3}V(\varepsilon z)u_n^2\psi_{\rho}\,{\rm d}z&+\int_{\mathbb{R}^3}\phi_{u_n}^tu_n^2\psi_{\rho}\,{\rm d}z\\
&=\int_{\mathbb{R}^3}g(\varepsilon z,u_n)u_n\psi_{\rho}\,{\rm d}z+o(1).
\end{align*}
Choose $\rho>0$ large enough such that $\Lambda'/\varepsilon\subset B_{\rho}(0)$, then $\varepsilon z\in\Lambda'$, by $(g_3)$, we have
\begin{align}\label{equ3-14}
&(1-\frac{1}{k})\int_{\mathbb{R}^3}V(\varepsilon z)u_n^2\psi_{\rho}\,{\rm d}z\nonumber\\
&\leq\int_{\mathbb{R}^3}|(-\Delta)^{\frac{s}{2}}u_n|^2\psi_{\rho}\,{\rm d}z+\int_{\mathbb{R}^3}V(\varepsilon z)u_n^2\psi_{\rho}\,{\rm d}z+\int_{\mathbb{R}^3}\phi_{u_n}^tu_n^2\psi_{\rho}\,{\rm d}z-\int_{\mathbb{R}^3}g(\varepsilon z,u_n)u_n\psi_{\rho}\,{\rm d}z\nonumber\\
&=\int_{\mathbb{R}^3}|(-\Delta)^{\frac{s}{2}}u_n|^2\psi_{\rho}\,{\rm d}z-\int_{\mathbb{R}^3}(-\Delta)^{\frac{s}{2}}u_n(-\Delta)^{\frac{s}{2}}(\psi_{\rho}u_n)\,{\rm d}z+o(1).
\end{align}
Now, by the nonlocal Leibniz rule, \eqref{equ2-5}, \eqref{equ2-6}, using H\"{o}lder's inequality and Proposition 3.4 in \cite{NPV}, we have that
\begin{align*}
&\Big|\int_{\mathbb{R}^3}(-\Delta)^{\frac{s}{2}}u_n(-\Delta)^{\frac{s}{2}}(\psi_{\rho}u_n)\,{\rm d}z-\int_{\mathbb{R}^3}|(-\Delta)^{\frac{s}{2}}u_n|^2\psi_{\rho}\,{\rm d}z\Big|\leq\frac{C}{\rho^s}\int_{\mathbb{R}^3}|u_n(-\Delta)^{\frac{s}{2}}u_n|\,{\rm d}z\\
&+\frac{C}{\rho^{\frac{s}{2}}}\int_{\mathbb{R}^3}|(-\Delta)^{\frac{s}{2}}u_n|\Big(\int_{\mathbb{R}^3}\frac{|u_n(z)-u_n(y)|^2}{|z-y|^{3+s}}\,{\rm d}y\Big)^{\frac{1}{2}}\,{\rm d}z\\
&\leq\frac{C}{\rho^s}\|u_n\|_{H_{\varepsilon}}^2+\frac{C}{\rho^{\frac{s}{2}}}\|u_n\|_{\mathcal{D}^{s,2}}\Big(\int_{\mathbb{R}^3}\int_{\mathbb{R}^3}\frac{|u_n(z)-u_n(y)|^2}{|z-y|^{3+s}}\,{\rm d}y\,{\rm d}z\Big)^{\frac{1}{2}}\\
&=\frac{C}{\rho^s}\|u_n\|_{H_{\varepsilon}}^2+\frac{C}{\rho^{\frac{s}{2}}}\|u_n\|_{\mathcal{D}^{s,2}}\Big(\int_{\mathbb{R}^3}|\xi|^s|\widehat{u}_n(\xi)|^2\,{\rm d}\xi\Big)^{\frac{1}{2}}\\
&\leq\frac{C}{\rho^s}\|u_n\|_{H_{\varepsilon}}^2+\frac{C}{\rho^{\frac{s}{2}}}\|u_n\|_{\mathcal{D}^{s,2}}\Big(\int_{\mathbb{R}^3}(1+|\xi|^{2s})|\widehat{u}_n(\xi)|^2\,{\rm d}\xi\Big)^{\frac{1}{2}}\\
&\leq C(\frac{1}{\rho^s}+\frac{1}{\rho^{\frac{s}{2}}})\|u_n\|_{H_{\varepsilon}}^2\leq\frac{C}{\rho^{\frac{s}{2}}}\|u_n\|_{H_{\varepsilon}}^2.
\end{align*}
In view of \eqref{equ3-14}, we have
\begin{align*}
(1-\frac{1}{k})\int_{\mathbb{R}^3}V(\varepsilon z)u_n^2\psi_{\rho}\,{\rm d}z\leq\frac{C}{\rho^{\frac{s}{2}}}\|u_n\|_{H_{\varepsilon}}^2+o(1)\leq\frac{C}{\rho^{\frac{s}{2}}}+o(1).
\end{align*}
Hence, we get
\begin{align}\label{equ3-15}
\limsup_{n\rightarrow\infty}\int_{|z|\geq2\rho}u_n^2\,{\rm d}z\leq\frac{C}{\rho^{\frac{s}{2}}}.
\end{align}
If $\{z_n\}$ is unbounded, by Lemma \ref{lem3-4} and \eqref{equ3-15}, we have
\begin{equation*}
0<\beta\leq\frac{C}{\rho^{\frac{s}{2}}}.
\end{equation*}
which achieves a contradiction for large $\rho$.
\end{proof}

\begin{proposition}\label{pro3-1}
The functional $J_{\varepsilon}$ possesses a nontrivial critical point $u_{\varepsilon}\in H_{\varepsilon}$ such that
\begin{equation*}
J_{\varepsilon}(u_{\varepsilon})=\inf_{\gamma\in \Gamma}\max_{t\in[0,1]}J_{\varepsilon}(\gamma(t))=\inf_{u\in H_{\varepsilon}\backslash\{0\}}\max_{t\geq0}J_{\varepsilon}(tu)=\inf_{u\in \mathcal{N}_{\varepsilon}}J_{\varepsilon}(u)
\end{equation*}
\end{proposition}

\begin{proof}
From Lemma \ref{lem3-3}, up to a subsequence, we may assume that there is $u:=u_{\varepsilon}\in H_{\varepsilon}$ such that $u_n\rightharpoonup u$ in $H_{\varepsilon}$, $u_n\rightarrow u$ in $L_{loc}^r(\mathbb{R}^3)$ for all $1\leq r<2_s^{\ast}$ and $u_n\rightarrow u$ a.e. in $\mathbb{R}^3$. Lemma \ref{lem3-4} and Lemma \ref{lem3-5} imply that $u$ is nontrivial. Moreover, by $(iv)$ of Lemma \ref{lem2-1}, it is easy to check that for any $\varphi\in H_{\varepsilon}$, $\langle J_{\varepsilon}'(u_n),\varphi\rangle\rightarrow\langle J_{\varepsilon}'(u),\varphi\rangle=0$ as $n\rightarrow\infty$, that is $u$ is a nontrivial critical point of $J_{\varepsilon}$. Next, we show that $J_{\varepsilon}(u)=c_{\varepsilon}$. Indeed, using the fact that $u\in\mathcal{N}_{\varepsilon}$, Fatou's Lemma, $(g_4)$ and \eqref{equ3-2}, we have
\begin{align}\label{equ3-16}
c_{\varepsilon}&\leq J_{\varepsilon}(u)=J_{\varepsilon}(u)-\frac{1}{4}\langle J_{\varepsilon}'(u),u\rangle\nonumber\\
&=\frac{1}{4}\|u\|_{H_{\varepsilon}}^2+\int_{\mathbb{R}^3}\Big(\frac{1}{4}g(\varepsilon z,u)u-G(\varepsilon z,u)\Big)\,{\rm d}z\nonumber\\
&=\frac{1}{4}\|u\|_{\mathcal{D}^{s,2}}^2+\frac{1}{4}\int_{\mathbb{R}^3}\frac{V(\varepsilon z)}{2}u^2\,{\rm d}z+\int_{\mathbb{R}^3}\Big(\frac{V(\varepsilon z)}{4}u^2+\frac{1}{4}g(\varepsilon z,u)u-G(\varepsilon z,u)\Big)\,{\rm d}z\nonumber\\
&\leq\liminf_{n\rightarrow\infty}\Big[\frac{1}{4}\|u_n\|_{\mathcal{D}^{s,2}}^2+\frac{1}{4}\int_{\mathbb{R}^3}\frac{V(\varepsilon z)}{2}u_n^2\,{\rm d}z+\int_{\mathbb{R}^3}\Big(\frac{V(\varepsilon z)}{4}u_n^2+\frac{1}{4}g(\varepsilon z,u_n)u_n\nonumber\\
&-G(\varepsilon z,u_n)\Big)\,{\rm d}z\Big]\nonumber\\
&=\liminf_{n\rightarrow\infty}\Big(J_{\varepsilon}(u_n)-\frac{1}{4}\langle J_{\varepsilon}'(u_n),u_n\rangle\Big)=c_{\varepsilon}.
\end{align}
The proof is completed.

\end{proof}

\begin{remark}\label{rem3-1}
From \eqref{equ3-16}, it is not difficult to verify that $\|u_n\|_{H_{\varepsilon}}\rightarrow\|u\|_{H_{\varepsilon}}$ as $n\rightarrow\infty$. Using the Brezis-Lieb Lemma, we conclude that $u_n\rightarrow u$ in $H_{\varepsilon}$.
\end{remark}

\section{The limiting problem}
In this section, we consider the limiting problem
\begin{equation}\label{equ4-1}
\left\{
  \begin{array}{ll}
   (-\Delta)^su+\mu u+\phi u=f(u)+u^{2_s^{\ast}-1} & \hbox{in $\mathbb{R}^3$,} \\
    (-\Delta)^t\phi=u^2, u>0& \hbox{in $\mathbb{R}^3$,}
  \end{array}
\right.
\end{equation}
where $\mu$ is a positive constant. The energy functional corresponding to problem \eqref{equ3-1} is
\begin{align*}
\mathcal{I}_{\mu}(u)&=\frac{1}{2}\int_{\mathbb{R}^3}|(-\Delta)^{\frac{s}{2}}u|^2\,{\rm d}x+\frac{\mu}{2}\int_{\mathbb{R}^3}|u|^2\,{\rm d}x+\frac{1}{4}\int_{\mathbb{R}^3}\phi_u^tu^2\,{\rm d}x-\int_{\mathbb{R}^3}F(u)\,{\rm d}x\\
&-\int_{\mathbb{R}^3}(u^{+})^{2_s^{\ast}}\,{\rm d}x,\quad u\in H^s(\mathbb{R}^3).
\end{align*}
The main result of this section is stated as follows.
\begin{proposition}\label{pro4-1}
Suppose that $f$ satisfies $(f_0)$-$(f_2)$, then autonomous problem \eqref{equ4-1} has a positive ground state solution $u\in H^s(\mathbb{R}^3)\cap C^{2,\alpha}(\mathbb{R}^3)$, such that $\mathcal{I}_{\mu}(u)=c_{\mu}$, where
\begin{equation*}
c_{\mu}=\inf_{\gamma\in\Gamma_{\mu}}\max_{\tau\in[0,1]}\mathcal{I}_{\mu}(\gamma(\tau))=\inf_{\mathcal{N}}\mathcal{I}_{\mu}(u)=\inf_{u\in H^s(\mathbb{R}^3)\backslash\{0\}}\max_{\tau\geq0}\mathcal{I}_{\mu}(\tau u),
\end{equation*}
where
\begin{equation*}
\mathcal{N}=\{u\in H^s(\mathbb{R}^3\backslash\{0\}\,\,\Big|\,\, \langle\mathcal{I}_{\mu}'(u),u\rangle=0\}
\end{equation*}
and
\begin{equation*}
\Gamma_{\mu}=\Big\{\gamma\in C([0,1],H^s(\mathbb{R}^3))\,\,\Big|\,\, \gamma(0)=0,\,\, \mathcal{I}_{\mu}(\gamma(1))<0\Big\}.
\end{equation*}
\end{proposition}

For proving Proposition \ref{pro4-1}, we will shoe the following preliminary results.
\begin{lemma}\label{lem4-1}
$\mathcal{I}_{\mu}$ possesses the mountain pass geometry:\\
$(i)$ there exist $\rho_0,\alpha_0>0$ such that $\mathcal{I}_{\mu}(u)\geq\alpha_0$ for all $u\in H^s(\mathbb{R}^3)$ with $\|w\|=\rho_0$;\\
$(ii)$ there exists $u_0\in H^s(\mathbb{R}^3)$ such that $\mathcal{I}_{\mu}(u_0)<0$.
\end{lemma}

Similar argument to Lemma \ref{lem3-2}, we can obtain the estimate of $c_{\mu}$.
\begin{lemma}\label{lem4-2}
Suppose $(f_0)-(f_2)$ hold, then $c_{\mu}$ satisfies
\begin{equation*}
0<c_{\mu}<\frac{s}{3}\mathcal{S}_s^{\frac{3}{2s}},
\end{equation*}
where $\mathcal{S}_s$ is the best Sobolev constant for the embedding $\mathcal{D}^{s,2}(\mathbb{R}^3)\hookrightarrow L^{2_s^{\ast}}(\mathbb{R}^3)$.
\end{lemma}

We can prove the following compactness condition.
\begin{lemma}\label{lem4-3}
Let $\{u_n\}$ be a $(PS)_c$ sequence for functional $\mathcal{I}_{\mu}$, then for any $c\in(0,\frac{s}{3}\mathcal{S}_s^{\frac{3}{2s}})$, under a translation, the sequence $\{u_n\}$ strongly convergence in $H^s(\mathbb{R}^3)$.
\end{lemma}
\begin{proof}
Suppose $\{u_n\}\subset H^s(\mathbb{R}^3)$ satisfies
\begin{equation}\label{equ4-2}
\mathcal{I}_{\mu}(u_n)\rightarrow c,\quad \mathcal{I}_{\mu}'(u_n)\rightarrow0\quad \text{as}\,\, n\rightarrow\infty.
\end{equation}
It follows from $(f_1)$-$(f_2)$ and \eqref{equ4-2} that
\begin{align*}
c+\|u_n\|&\geq\mathcal{I}_{\mu}(u_n)-\frac{1}{4}\langle\mathcal{I}_{\mu}'(u_n),u_n\rangle\\
&=\frac{1}{4}\|u_n\|^2+\int_{\mathbb{R}^3}\Big(\frac{1}{4}f(u_n)u_n-F(u_n)\Big)\,{\rm d}x+\frac{4s-3}{12}\int_{\mathbb{R}^3}(u_n^{+})^{2_s^{\ast}}\,{\rm d}x\\
&\geq\frac{1}{4}\|u_n\|^2.
\end{align*}
Thus $\{u_n\}$ is bounded in $H^s(\mathbb{R}^3)$. Similar argument to the proof of Lemma \ref{lem3-4}, we can obtain the bounded sequence $\{u_n\}$ is nonvanishing, i.e., there exist a sequence $\{x_n\}\subset\mathbb{R}^3$ and $R_0>0$, $\beta_0>0$ such that
\begin{equation}\label{equ4-3}
\liminf_{n\rightarrow\infty}\int_{B_{R_0}(x_n)}|u_n|^2\,{\rm d}z\geq \beta_0.
\end{equation}
Up to a subsequence, we may assume that there is $u\in H^s(\mathbb{R}^3)$ such that $u_n\rightharpoonup u$ in $H^s(\mathbb{R}^3)$. By \eqref{equ4-3}, we can assume that $u\neq0$. Indeed, if $u=0$, then $u_n\rightharpoonup 0$ in $H^s(\mathbb{R}^3)$ and $u_n\rightarrow0$ (otherwise, contradicts with $c>0$). Set $v_n(x)=u_n(x+x_n)$, then $\{v_n\}$ is also a $(PS)_c$ sequence for $\mathcal{I}_{\mu}$, so $\{v_n\}$ is bounded in $H^s(\mathbb{R}^3)$ and there exists $v\in H^s(\mathbb{R}^3)$ such that $v_n\rightharpoonup v$ in $H^s(\mathbb{R}^3)$, $v_n\rightarrow v$ in $L_{loc}^r(\mathbb{R}^3)$ for all $1\leq r<2_s^{\ast}$ and $v_n\rightarrow v$ a.e. in $\mathbb{R}^3$. By $(iv)$ of Lemma \ref{lem2-1}, we have that $\langle\mathcal{I}_{\mu}'(v_n),\varphi\rangle\rightarrow\langle\mathcal{I}_{\mu}'(v),\varphi\rangle=0$ for any $\varphi\in H^s(\mathbb{R}^3)$. Clearly, $\langle\mathcal{I}_{\mu}'(v),v\rangle=0$.

Next we only need to show that $v_n\rightarrow v$ in $H^s(\mathbb{R}^3)$. By the weakly semi-lower continuity of norm, we have that
\begin{equation}\label{equ4-4}
\liminf\limits_{n\rightarrow\infty}\|v_n\|\geq\|v\|.
\end{equation}
In order to prove that \eqref{equ4-4} hold, we must show the equality holds in \eqref{equ4-4}. Otherwise, by Fatou's Lemma, we get
\begin{align*}
c\leq\mathcal{I}_{\mu}(v)&=\mathcal{I}_{\mu}(v)-\frac{1}{4}\langle\mathcal{I}_{\mu}'(v),v\rangle\\
&=\frac{1}{4}\|v\|^2+\int_{\mathbb{R}^3}(\frac{1}{4}f(v)v-F(v))\,{\rm d}x+\frac{4s-3}{12}\int_{\mathbb{R}^3}(v^{+})^{2_s^{\ast}}\,{\rm d}x\\
&<\liminf_{n\rightarrow\infty}\Big[\frac{1}{4}\|v\|^2+\int_{\mathbb{R}^3}(\frac{1}{4}f(v)v-F(v))\,{\rm d}x+\frac{4s-3}{12}\int_{\mathbb{R}^3}(v^{+})^{2_s^{\ast}}\,{\rm d}x\Big]\\
&=\liminf_{n\rightarrow\infty}\Big(\mathcal{I}_{\mu}(v_n)-\frac{1}{4}\langle\mathcal{I}_{\mu}'(v_n),v_n\rangle\Big)\\
&=c
\end{align*}
which is a contradiction. Thus, up to a subsequence, using Brezis-Lieb Lemma, we conclude that $v_n\rightarrow v$ in $H^s(\mathbb{R}^3)$.

{\bf Proof of Proposition \ref{pro4-1}.} From Lemma \ref{lem4-1}, and using the mountain-pass Lemma without $(PS)$ condition, we get a $(PS)_{c_{\mu}}$ sequence $\{u_n\}\subset H^s(\mathbb{R}^3)$. By Lemma \ref{lem4-2} and Lemma \ref{lem4-3}, under a translation, still denoted by $\{u_n\}$, there is $u\in H^s(\mathbb{R}^3)$ such that $u_n\rightarrow u$ in $H^s(\mathbb{R}^3)$. Therefore, $\mathcal{I}_{\mu}(u)=c_{\mu}$ and $\mathcal{I}_{\mu}'(u)=0$, i.e., $u\in H^s(\mathbb{R}^3)$ is a weak solution of problem \eqref{equ4-1}. By standard argument to the proof Proposition 4.4 in \cite{Teng2}, we have that $u\in C^{2,\alpha}(\mathbb{R}^3)$ for some $\alpha\in(0,1)$. The remain proof is to show $u$ is positive. Using $-u^{-}$ as a testing function, it is easy to see that $u\geq0$. Since $u\in C^{2,\alpha}(\mathbb{R}^3)$, by Lemma 3.2 in \cite{NPV}, we have that
\begin{equation*}
(-\Delta)^su(x)=-\frac{1}{2}C(s)\int_{\mathbb{R}^3}\frac{u(x+y)+u(x-y)-2u(x)}{|x-y|^{3+2s}}\,{\rm d}x\,{\rm d}y,\quad \forall\,\, x\in\mathbb{R}^3.
\end{equation*}
Assume that there exists $x_0\in\mathbb{R}^3$ such that $u(x_0)=0$, then from $u\geq0$ and $u\not\equiv0$, we get
\begin{equation*}
(-\Delta)^su(x_0)=-\frac{1}{2}C(s)\int_{\mathbb{R}^3}\frac{u(x_0+y)+u(x_0-y)}{|x_0-y|^{3+2s}}\,{\rm d}x\,{\rm d}y<0.
\end{equation*}
However, observe that $(-\Delta)^su(x_0)=-\mu u(x_0)-(\phi_u^tu)(x_0)+f(u(x_0))+u(x_0)^{2_s^{\ast}-1}=0$, a contradiction. Hence, $u(x)>0$, for every $x\in\mathbb{R}^3$. The proof is completed.
\end{proof}

For $V_0=\min\limits_{\Lambda}V$, let $w$ be a ground state solution to the following problem
\begin{equation}\label{equ4-5}
\left\{
  \begin{array}{ll}
   (-\Delta)^su+V_0 u+\phi u=f(u)+u^{2_s^{\ast}-1} & \hbox{in $\mathbb{R}^3$,} \\
    (-\Delta)^t\phi=u^2, u>0& \hbox{in $\mathbb{R}^3$,}
  \end{array}
\right.
\end{equation}
satisfying $\mathcal{I}_{V_0}(w)=\inf\limits_{v\in H^s(\mathbb{R}^3)\backslash\{0\}}\max\limits_{\tau\geq0}\mathcal{I}_{V_0}(\tau v):=c_{V_0}$.

\begin{lemma}\label{lem4-4}
\begin{equation}\label{equ4-6}
c_{\varepsilon}\leq c_{V_0}+o(1).
\end{equation}
\end{lemma}

\begin{proof}
Let $z_0\in\Lambda$ be such that $V(z_0)=V_0$ and so there is $R_0>0$ such that$B_{2R_0}(z_0)\subset \Lambda$. Let $v_{\varepsilon}(z)=\eta(\varepsilon z-z_0)w(\frac{\varepsilon z-z_0}{\varepsilon})$, where $\eta$ is smooth cut-off function satisfying $0\leq\eta\leq1$, $\eta=1$ on $B_{R_0}(0)$, $\eta=0$ on $\mathbb{R}^3\backslash B_{2R_0}(0)$, $|\nabla\eta|\leq C$. By using the change of variables $\varepsilon z=\varepsilon z'-z_0$, we can write
\begin{align*}
&J_{\varepsilon}(tv_{\varepsilon})=\frac{t^2}{2}\int_{\mathbb{R}^3}|(-\Delta)^{\frac{s}{2}}v_{\varepsilon}|^2\,{\rm d}z+\frac{t^2}{2}\int_{\mathbb{R}^3}V(\varepsilon z+z_0)\eta(\varepsilon z)w(z)\,{\rm d}z\\
&+\frac{t^4}{4}\int_{\mathbb{R}^3}\phi_{\eta(\varepsilon\cdot) w}^t\eta(\varepsilon z)w(z)\,{\rm d}z-\int_{\mathbb{R}^3}F(t\eta(\varepsilon z)w(z))\,{\rm d}z-\frac{t^{2_s^{\ast}}}{2_s^{\ast}}\int_{\mathbb{R}^3}(\eta(\varepsilon z)w^{+}(z))^{2_s^{\ast}}\,{\rm d}z
\end{align*}

Since $w>0$, by standard argument, there is a unique $t_{\varepsilon}>0$ such that $\sup\limits_{t>0}J_{\varepsilon}(tv_{\varepsilon})=J_{\varepsilon}(t_{\varepsilon}v_{\varepsilon})$ and $\frac{d J_{\varepsilon}(tv_{\varepsilon})}{dt}|_{t=t_{\varepsilon}}=0$, i.e.,
\begin{align*}
t_{\varepsilon}\int_{\mathbb{R}^3}(|(-\Delta)^{\frac{s}{2}}v_{\varepsilon}|^2&+V(\varepsilon z+z_0)\eta(\varepsilon z) w^2(z))\,{\rm d}z+t_{\varepsilon}^3\int_{\mathbb{R}^3}\phi_{\eta(\varepsilon\cdot) w}^t\eta(\varepsilon z) w^2(z)\,{\rm d}z\\
&=\int_{\mathbb{R}^3} f(t_{\varepsilon}\eta(\varepsilon z) w(z))\eta(\varepsilon z) w(z),{\rm d}z+t_{\varepsilon}^{2_s^{\ast}-1}\int_{\mathbb{R}^3}(\eta(\varepsilon z) w^{+}(z))^{2_s^{\ast}}\,{\rm d}z
\end{align*}
which means that $t_{\varepsilon}v_{\varepsilon}\in\mathcal{N}_{\varepsilon}$. We claim that there exist $T_1,T_2>0$ such that $0<T_1\leq t_{\varepsilon}\leq T_2$.
Observe that
\begin{equation}\label{equ4-7}
\lim_{\varepsilon\rightarrow0}\int_{\mathbb{R}^3}(|(-\Delta)^{\frac{s}{2}}v_{\varepsilon}|^2+V(\varepsilon z+z_0)\eta(\varepsilon z) w^2(z))\,{\rm d}z=\int_{\mathbb{R}^3}(|(-\Delta)^{\frac{s}{2}}w|^2+V(z_0)w^2)\,{\rm d}z,
\end{equation}
\begin{equation}\label{equ4-8}
\lim_{\varepsilon\rightarrow0}\int_{\mathbb{R}^3}\phi_{\eta(\varepsilon\cdot) w}^t\eta(\varepsilon z) w^2(z)\,{\rm d}z=\int_{\mathbb{R}^3}\phi_{w}^t w^2\,{\rm d}z
\end{equation}
and
\begin{equation}\label{equ4-9}
\lim_{\varepsilon\rightarrow0}\int_{\mathbb{R}^3}|\eta(\varepsilon z) w(z)|^r\,{\rm d}z=\int_{\mathbb{R}^3}|w|^r\,{\rm d}z\quad \text{for all}\,\,2\leq r\leq2_s^{\ast}.
\end{equation}
If $t_{\varepsilon}\rightarrow0$ as $\varepsilon\rightarrow0$, by $(f_0)$ and $(f_2)$, we have that
\begin{align*}
\int_{\mathbb{R}^3}(|(-\Delta)^{\frac{s}{2}}v_{\varepsilon}|^2&+V(\varepsilon z+z_0)\eta(\varepsilon z) w^2(z))\,{\rm d}z+t_{\varepsilon}^2\int_{\mathbb{R}^3}\phi_{\eta(\varepsilon\cdot) w}^t\eta(\varepsilon z) w^2\,{\rm d}z\\
&\leq Ct_{\varepsilon}^{p-2}\int_{\mathbb{R}^3}|\eta(\varepsilon z) w(z)|^p\,{\rm d}z+t_{\varepsilon}^{2_s^{\ast}-2}\int_{\mathbb{R}^3}|\eta(\varepsilon z) w(z)|^{2_s^{\ast}}\,{\rm d}z
\end{align*}
which leads to a contradiction using \eqref{equ4-7}-\eqref{equ4-9}. If $t_{\varepsilon}\rightarrow\infty$ as $\varepsilon\rightarrow0$, then
\begin{align*}
\int_{\mathbb{R}^3}(|(-\Delta)^{\frac{s}{2}}v_{\varepsilon}|^2&+V(\varepsilon z+z_0)\eta(\varepsilon z) w^2(z))\,{\rm d}z+t_{\varepsilon}^2\int_{\mathbb{R}^3}\phi_{\eta(\varepsilon\cdot) w}^t\eta(\varepsilon z) w^2(z)\,{\rm d}z\\
&\geq t_{\varepsilon}^{2_s^{\ast}-2}\int_{\mathbb{R}^3}(\eta(\varepsilon z) w^{+}(z))^{2_s^{\ast}}\,{\rm d}z.
\end{align*}
Hence, by \eqref{equ4-7}-\eqref{equ4-9}, it is easy to achieve a contradiction. Thus the claim holds. Going if necessary to a subsequence, we may assume that $t_{\varepsilon}\rightarrow T>0$ as $\varepsilon\rightarrow0$, then using \eqref{equ4-7}-\eqref{equ4-9}, we get
\begin{align*}
T\int_{\mathbb{R}^3}(|(-\Delta)^{\frac{s}{2}}w|^2+V_0w^2)\,{\rm d}z+T^3\int_{\mathbb{R}^3}\phi_{w}^tw^2\,{\rm d}z&=\int_{\mathbb{R}^3} f(Tw)w\,{\rm d}z\\
&+T^{2_s^{\ast}-1}\int_{\mathbb{R}^3}(w^{+})^{2_s^{\ast}}\,{\rm d}z.
\end{align*}
Since $w$ is a solution of problem \eqref{equ4-1}, we have that
\begin{align*}
(\frac{1}{T^2}-1)\int_{\mathbb{R}^3}(|(-\Delta)^{\frac{s}{2}}w|^2+V_0w^2)\,{\rm d}z&=\int_{\mathbb{R}^3} w^4\Big(\frac{f(Tw)}{(Tw)^3}-\frac{f(w)}{w^3}\Big)\,{\rm d}z\\
&+(T^{2_s^{\ast}-4}-1)\int_{\mathbb{R}^3}(w^{+})^{2_s^{\ast}}\,{\rm d}z.
\end{align*}
By $(f_2)$, we see that $t_{\varepsilon}\rightarrow T=1$. Therefore, by \eqref{equ4-7}-\eqref{equ4-9}, we get
\begin{align*}
c_{\varepsilon}\leq J_{\varepsilon}(t_{\varepsilon}v_{\varepsilon})=\mathcal{I}_{V_0}(w)+o(1)=c_{V_0}+o(1).
\end{align*}
Thus \eqref{equ4-6} follows.
\end{proof}

\section{Uniformly estimate of solution sequence.}
In this section, we consider the following problem
\begin{equation}\label{equ5-1}
\left\{
  \begin{array}{ll}
   (-\Delta)^su+V_n(x) u+ \phi u=f_n(x,u) & \hbox{in $\mathbb{R}^3$,} \\
    (-\Delta)^t\phi=u^2& \hbox{in $\mathbb{R}^3$,}
  \end{array}
\right.
\end{equation}
where $\{V_n\}$ satisfies $V_n(x)\geq \alpha_0>0$ for all $x\in\mathbb{R}^3$ and $f_n(x,\tau)$ is a Carathedory function satisfying that for any $\delta>0$, there exists $C_{\varepsilon}>0$ such that
\begin{equation}\label{equ5-2}
|f_n(x,\tau)|\leq\delta|\tau|+C_{\delta}|\tau|^{2_s^{\ast}-1},\quad \forall (x,\tau)\in\mathbb{R}^3\times\mathbb{R}.
\end{equation}

\begin{proposition}\label{pro5-1}
Assume that $u_n$ are nonnegative weak solution of \eqref{equ5-1} satisfying $u_n$ convergence strongly in $H^s(\mathbb{R}^3)$. Then there exists $C>0$ such that
\begin{equation*}
\|u_n\|_{L^{\infty}}\leq C\quad \text{for all}\,\, n.
\end{equation*}
\end{proposition}
\begin{proof}
Define
\begin{equation*}
\psi_{T,k}(t)=\left\{
          \begin{array}{ll}
            0 & \hbox{if $t\leq0$,} \\
            t^k & \hbox{if $0<t<T$,} \\
            kT^{k-1}(t-T)+T^k, & \hbox{if $t\geq T$.}
          \end{array}
        \right.
\end{equation*}
Clearly $\psi$ is a convex and differentiable function, and $(-\Delta)^s\psi_{T,k}(u)\leq\psi_{T,k}'(u)(-\Delta)^su$. Moreover, $\|\psi_{T,k}(u)\|_{\mathcal{D}^{s,2}}\leq kT^{k-1}\|u\|_{\mathcal{D}^{s,2}}$. Taking $\psi_{T,k}(u_n)\psi_{T,k}'(u_n)$ as a test function in \eqref{equ5-1}, by \eqref{equ5-2}, we get
\begin{align*}
&\int_{\mathbb{R}^3}(-\Delta)^su_n\psi_{T,k}(u_n)\psi_{T,k}'(u_n)\,{\rm d}x+\int_{\mathbb{R}^3}V_n(x)u_n\psi_{T,k}(u_n)\psi_{T,k}'(u_n)\,{\rm d}x\\
&+\int_{\mathbb{R}^3}\phi_{u_n}^tu_n\psi_{T,k}(u_n)\psi_{T,k}'(u_n)\,{\rm d}x=\int_{\mathbb{R}^3}f_n(x,u_n)\psi_{T,k}(u_n)\psi_{T,k}'(u_n)\,{\rm d}x\\
&\leq\delta\int_{\mathbb{R}^3}u_n\psi_{T,k}(u_n)\psi_{T,k}'(u_n)\,{\rm d}x+C_{\delta}\int_{\mathbb{R}^3}u_n^{2_s^{\ast}-1}\psi_{T,k}(u_n)\psi_{T,k}'(u_n)\,{\rm d}x
\end{align*}
Taking $\delta=\alpha_0$, using Sobolev inequality, H\"{o}lder's inequality and the fact $t\psi_{T,k}'(t)\leq k\psi_{T,k}(t)$ for $t\geq0$, we have that
\begin{align}\label{equ5-3}
&\|\psi_{T,k}(u_n)\|_{L^{2_s^{\ast}}}^2\leq\mathcal{S}_s^{-1}\|\psi_{T,k}(u_n)\|_{\mathcal{D}^{s,2}}^2=\mathcal{S}_s^{-1}\int_{\mathbb{R}^3}\psi_{T,k}(u_n)(-\Delta)^s\psi_{T,k}(u_n)\,{\rm d}x\nonumber\\
&\leq\mathcal{S}_s^{-1}\int_{\mathbb{R}^3}(-\Delta)^su_n\psi_{T,k}(u_n)\psi_{T,k}'(u_n)\,{\rm d}x\nonumber\\
&\leq C\int_{\mathbb{R}^3}u_n^{2_s^{\ast}-1}\psi_{T,k}(u_n)\psi_{T,k}'(u_n)\,{\rm d}x\leq Ck\int_{\mathbb{R}^3}u_n^{2_s^{\ast}-2}\psi_{T,k}^2(u_n)\,{\rm d}x\\
&\leq Ck\Big[\int_{u_n< A_0}u_n^{2_s^{\ast}-2}\psi_{T,k}^2(u_n)\,{\rm d}x+\|\psi_{T,k}(u_n)\|_{L^{2_s^{\ast}}}^2\Big(\int_{u_n\geq A_0}u_n^{2_s^{\ast}}\,{\rm d}x\Big)^{\frac{2_s^{\ast}-2}{2_s^{\ast}}}\Big].\nonumber
\end{align}
Take $k=k_1=\frac{2_s^{\ast}}{2}$, we get
\begin{equation}\label{equ5-4}
\|\psi_{T,k_1}(u_n)\|_{L^{2_s^{\ast}}}^2\leq C\frac{2_s^{\ast}}{2}\Big[\int_{u_n< A_0}u_n^{2_s^{\ast}-2}\psi_{T,k_1}^2(u_n)\,{\rm d}x+\|\psi_{T,k_1}(u_n)\|_{L^{2_s^{\ast}}}^2\Big(\int_{u_n\geq A_0}u_n^{2_s^{\ast}}\,{\rm d}x\Big)^{\frac{2_s^{\ast}-2}{2_s^{\ast}}}\Big].
\end{equation}
Since $u_n$ convergence strongly in $H^s(\mathbb{R}^3)$, then $u_n$ convergence strongly in $L^{2_s^{\ast}}(\mathbb{R}^3)$. Thus, we can take $A_0$ large enough such that
\begin{equation*}
\Big(\int_{u_n\geq A_0}u_n^{2_s^{\ast}}\,{\rm d}x\Big)^{\frac{2_s^{\ast}-2}{2_s^{\ast}}}\leq\frac{1}{C2_s^{\ast}}.
\end{equation*}
Thus, from $\psi_{T,k}(t)\leq t^k$ for $t\geq0$ and \eqref{equ5-4}, we deduce that
\begin{equation}\label{equ5-5}
\|\psi_{T,k_1}(u_n)\|_{L^{2_s^{\ast}}}^2\leq CA_0^{2_s^{\ast}-2}\int_{u_n< A_0}|u_n|^{2_s^{\ast}}\,{\rm d}x\leq CA_0^{2_s^{\ast}-2}\int_{\mathbb{R}^3}|u_n|^{2_s^{\ast}}\,{\rm d}x.
\end{equation}
Letting $T\rightarrow+\infty$ in \eqref{equ5-5}, we get
\begin{equation}\label{equ5-6}
\|u_n\|_{L^{2_s^{\ast}k_1}}^{2k_1}\leq CA_0^{2_s^{\ast}-2}\|u_n\|_{L^{2_s^{\ast}}}^{2_s^{\ast}}<\infty.
\end{equation}
Letting $T\rightarrow+\infty$ in \eqref{equ5-3}, we have
\begin{equation*}
\Big(\int_{\mathbb{R}^3}|u_n|^{2_s^{\ast}k}\,{\rm d}x\Big)^{\frac{1}{2^{\ast}(k-1)}}\leq (Ck)^{\frac{1}{2(k-1)}}\Big(\int_{\mathbb{R}^3}|u_n|^{2k+2_s^{\ast}-2}\,{\rm d}x\Big)^{\frac{1}{2(k-1)}}.
\end{equation*}
For $m\geq1$, we define $k_{m+1}$ inductively so that $2k_{m+1}+2_s^{\ast}-2=2_s^{\ast}k_m$ and $k_1=\frac{2_s^{\ast}}{2}$, using \eqref{equ5-6}, it is easy to check that
\begin{align*}
\Big(\int_{\mathbb{R}^3}|u_n|^{2_s^{\ast}k_{m+1}}\,{\rm d}x\Big)^{\frac{1}{2_s^{\ast}(k_{m+1}-1)}}&\leq\prod_{i=1}^m(Ck_{i+1})^{\frac{1}{2k_{i+1}-1}}\Big(\int_{\mathbb{R}^3}|u_n|^{2_s^{\ast}k_1}\,{\rm d}x\Big)^{\frac{1}{2_s^{\ast}(k_1-1)}}\\
&\leq C\prod_{i=1}^m(Ck_{i+1})^{\frac{1}{2k_{i+1}-1}}\|u_n\|_{L^{2_s^{\ast}}}^{\frac{(2_s^{\ast})^3(2_s^{\ast}-2)}{4}},
\end{align*}
letting $m\rightarrow\infty$ in the above inequality, we conclude that
\begin{equation*}
\|u_n\|_{L^{\infty}}\leq C\|u_n\|_{L^{2_s^{\ast}}}^{\frac{(2_s^{\ast})^3(2_s^{\ast}-2)}{4}}\leq C.
\end{equation*}

\end{proof}

\section{Proof of Theorem \ref{thm1-1}.}

For $\varepsilon>0$, let $u_{\varepsilon}$ be the mountain-pass solution to \eqref{equ3-1} given by Proposition \ref{pro3-1}. For any sequence $\{\varepsilon_n\}$ satisfying $\varepsilon_{n}\rightarrow0^{+}$, denote by $u_n:=u_{\varepsilon_n}$, $J_n:=J_{\varepsilon_n}$ and $H_{n}:=H_{\varepsilon_n}$. Then $u_n$ satisfies
\begin{equation}\label{equ6-1}
(-\Delta)^su_n+V(\varepsilon_n z) u_n+ \phi_{u_n}^t u_n=g(\varepsilon_n z,u_n)\quad z\in\mathbb{R}^3.
\end{equation}
Here $u_n$ is a critical point of $J_n$ and $J_n(u_n)=c_{\varepsilon_n}$. Using Lemma \ref{lem4-4}, and similar argument to the proof of Lemma \ref{lem3-3}, we have that $\{u_n\}$ is bounded in $H_n$. Similar to Lemma \ref{lem3-4}, we have
\begin{lemma}\label{lem6-1}
There exist a sequence $\{y_n\}\subset\mathbb{R}^3$ and $R>0$, $\beta>0$ such that
\begin{equation*}
\int_{B_R(y_n)}u_n^2\,{\rm d}z\geq\beta.
\end{equation*}
\end{lemma}

\begin{lemma}\label{lem6-2}
$\{\varepsilon_ny_n\}$ is bounded in $\mathbb{R}^3$. Moreover, ${\rm dist}(\varepsilon_ny_n,\Lambda')\leq\varepsilon_n R$.
\end{lemma}
\begin{proof}
For $\delta>0$, define $K_{\delta}=\{z\in\mathbb{R}^3\,\,|\,\,{\rm dist}(z,\Lambda')\leq\delta\}$. Let $\varphi_{\varepsilon_n}(z)=\varphi(\varepsilon_n z)$ satisfy that $\varphi\in C^{\infty}(\mathbb{R}^3,[0,1])$, $\varphi=1$ on $\mathbb{R}^3\backslash K_{\delta}$, $\varphi=0$ on $\Lambda'$ and $|\nabla\varphi|\leq\frac{C}{\delta}$. Noting that ${\rm supp}\varphi_{\varepsilon_n}\cap(\Lambda'/\varepsilon_n)=\emptyset$, then $g(\varepsilon_n z,u_n)u_n\varphi_{\varepsilon_n}=\tilde{f}(u_n)u_n\varphi_{\varepsilon_n}$. Taking $u_n\varphi_{\varepsilon_n}$ as a test function in \eqref{equ6-1}, by $(g_3)$, similar argument to Lemma \ref{lem3-5}, we have that
\begin{align*}
&\alpha(1-\frac{1}{k})\int_{\mathbb{R}^3}u_n^2\varphi_{\varepsilon_n}\,{\rm d}z\leq(1-\frac{1}{k})\int_{\mathbb{R}^3}V(\varepsilon_n z)u_n^2\varphi_{\varepsilon_n}\,{\rm d}z+\int_{\mathbb{R}^3}|(-\Delta)^{\frac{s}{2}}u_n|^2\varphi_{\varepsilon_n}\,{\rm d}x\\
&\leq(1-\frac{1}{k})\int_{\mathbb{R}^3}V(\varepsilon_n z)u_n^2\varphi_{\varepsilon_n}\,{\rm d}z+\int_{\mathbb{R}^3}\phi_{u_n}^tu_n^2\varphi_{\varepsilon_n}\,{\rm d}z+\int_{\mathbb{R}^3}|(-\Delta)^{\frac{s}{2}}u_n|^2\varphi_{\varepsilon_n}\,{\rm d}x\\
&\leq\int_{\mathbb{R}^3}|(-\Delta)^{\frac{s}{2}}u_n|^2\varphi_{\varepsilon_n}\,{\rm d}x-\int_{\mathbb{R}^3}(-\Delta)^{\frac{s}{2}}u_n(-\Delta)^{\frac{s}{2}}(u_n\varphi_{\varepsilon_n})\,{\rm d}z\leq C\frac{\varepsilon_n}{\delta^{s}}.
\end{align*}
We claim that for small $\varepsilon_n>0$, there exists $y_n'$ such that $\varepsilon_ny_n'\in K_{\delta}$ and $|y_n'-y_n|\leq R$. Otherwise, there exists a subsequence $\varepsilon_{n_j}\rightarrow0$ such that $|z-y_{n_j}|\leq R$ and $\varepsilon_{n_j}z\not\in K_{\delta}$, that is, $B_R(y_{n_j})\cap\{z\in\mathbb{R}^3\,|\,\varepsilon_{n_j}z\in K_{\delta}\}={\O}$. Thus,
\begin{equation*}
\alpha(1-\frac{1}{k})\int_{B_R(y_{n_j})}u_n^2\,{\rm d}z\leq C\frac{\varepsilon_n}{\delta^{s}}
\end{equation*}
which contradicts with Lemma \ref{lem6-1}. Thus the claim follows. Moreover,
\begin{equation*}
{\rm dist}(\varepsilon_ny_n,\Lambda')\leq|\varepsilon_ny_n-\varepsilon_ny_n'|+{\rm dist}(\varepsilon_ny_n',\Lambda')\leq\varepsilon_n R+\delta.
\end{equation*}
By the arbitrariness of $\delta$, we complete the proof.
\end{proof}

By Lemma \ref{lem6-2}, we see that $\lim\limits_{n\rightarrow\infty}{\rm dist}(\varepsilon_ny_n,\Lambda')=0$, hence, there is a subsequence of $\{\varepsilon_ny_n\}$, still denoted by $\varepsilon_ny_n$ and $x_0\in\overline{\Lambda'}$ such that $\lim\limits_{n\rightarrow\infty}\varepsilon_ny_n=x_0$. Set
\begin{align*}
\widetilde{J}_{V(x_0)}(u)=\frac{1}{2}\int_{\mathbb{R}^3}|(-\Delta)^{\frac{s}{2}}u|^2\,{\rm d}x+\frac{1}{2}\int_{\mathbb{R}^3}V(x_0)u^2\,{\rm d}x+\frac{1}{4}\int_{\mathbb{R}^3}\phi_{u}^tu^2\,{\rm d}x-\int_{\mathbb{R}^3}G(u)\,{\rm d}x
\end{align*}
where $G(u)=\int_0^ug(s)\,{\rm d}s$ and $g(u)=\chi(x_0)(f(u)+(u^{+})^{2_s^{\ast}-1})+(1-\chi(x_0))\tilde{f}(u)$.

\begin{lemma}\label{lem6-3}
$x_0\in\Lambda$.
\end{lemma}
\begin{proof}
It suffices to show that $V(x_0)=V_0$. If this fact is proved, by $(V_1)$ and the definition of $\Lambda'$, we see that $x_0\not\in\partial\Lambda$ and $x_0\not\in\overline{\Lambda'}/\Lambda$, then $x_0\in\Lambda$.

Now, clearly, $V(x_0)\geq V_0$. The remain is to prove $V(x_0)\leq V_0$. Set $v_n(z)=u_n(z+y_n)$, then $v_n$ satisfies
\begin{equation}\label{equ6-2}
(-\Delta)^sv_n+V(\varepsilon_n z+\varepsilon_ny_n) v_n+ \phi_{v_n}^t v_n=g(\varepsilon_n z+\varepsilon_ny_n,v_n)\quad z\in\mathbb{R}^3
\end{equation}
\begin{equation}\label{equ6-3}
\int_{B_R(0)}v_n^2\,{\rm d}z\geq\beta>0
\end{equation}
and $\|v_n\|=\|u_n\|$ is bounded. Up to a subsequence, there exists $v\in H^s(\mathbb{R}^3)\backslash\{0\}$ such that $v_n\rightharpoonup v$ in $H^s(\mathbb{R}^3)$, $v_n\rightarrow v$ in $L_{loc}^r(\mathbb{R}^3)$ for all $1\leq r<2_s^{\ast}$ and $v_n\rightarrow v$ a.e. in $\mathbb{R}^3$.

Therefore, by $(iv)$ of Lemma \ref{lem2-1}, it is easy to show that
\begin{align*}
\int_{\mathbb{R}^3}(-\Delta)^{\frac{s}{2}}v(-\Delta)^{\frac{s}{2}}\varphi\,{\rm d}x+\int_{\mathbb{R}^3}V(x_0)v\varphi\,{\rm d}x+\int_{\mathbb{R}^3}\phi_{v}^tv\varphi\,{\rm d}x
=\int_{\mathbb{R}^3}g(v)\varphi\,{\rm d}x
\end{align*}
for any $\varphi\in H^s(\mathbb{R}^3)$. Thus, we get that $\langle \widetilde{J}_{V(x_0)}'(v),v\rangle=0$.

Let $c_{x_0}$ be the mountain-pass energy of $\widetilde{J}_{V(x_0)}$. Since $\widetilde{J}_{V(x_0)}(v)\geq \mathcal{I}_{V(x_0)}(v)$, then $c_{x_0}\geq c_{V(x_0)}$. Thus, similar argument to \eqref{equ3-16}, we get
\begin{align}\label{equ6-4}
c_{V(x_0)}\leq c_{x_0}&\leq\widetilde{J}_{V(x_0)}(v)-\frac{1}{4}\langle \widetilde{J}_{V(x_0)}'(v),v\rangle\nonumber\\
&=\frac{1}{4}\int_{\mathbb{R}^3}(|(-\Delta)^{\frac{s}{2}}v|^2+V(x_0)v^2)\,{\rm d}z+\int_{\mathbb{R}^3}(\frac{1}{4}g(v)v-G(v))\,{\rm d}z\nonumber\\
&\leq\liminf_{n\rightarrow\infty}\Big[\frac{1}{4}\int_{\mathbb{R}^3}(|(-\Delta)^{\frac{s}{2}}v_n|^2+V(\varepsilon_n z+\varepsilon_ny_n)v_n^2)\,{\rm d}z\nonumber\\
&+\int_{\mathbb{R}^3}\Big(\frac{1}{4}g(\varepsilon_n z+\varepsilon_ny_n,v_n)v_n-G(\varepsilon_n z+\varepsilon_ny_n,v_n)\Big)\,{\rm d}z\Big]\nonumber\\
&=\liminf_{n\rightarrow\infty}\Big[\frac{1}{4}\int_{\mathbb{R}^3}(|(-\Delta)^{\frac{s}{2}}u_n|^2+V(\varepsilon_n z+\varepsilon_ny_n)u_n^2)\,{\rm d}z\nonumber\\
&+\int_{\mathbb{R}^3}\Big(\frac{1}{4}g(\varepsilon_n z,u_n)u_n-G(\varepsilon_n z,u_n)\Big)\,{\rm d}z\Big]\nonumber\\
&=\liminf_{n\rightarrow\infty}\Big(J_n(u_n)-\frac{1}{4}\langle J_n'(u_n),u_n\rangle\Big)\leq c_{V_0}.
\end{align}
Assume by the contrary that $V(x_0)>V_0$. Denote $w$ be a ground state critical point of $\mathcal{I}_{V(x_0)}$. By standard argument, there exists $\tau_0>0$ such that $\mathcal{I}_{V_0}(\tau_0w)=\sup_{\tau>0}\mathcal{I}_{V_0}(\tau w)$. Hence,
\begin{align*}
c_{V_0}\leq\sup_{\tau>0}\mathcal{I}_{V_0}(\tau w)=\mathcal{I}_{V_0}(\tau_0w)<\mathcal{I}_{V(x_0)}(\tau_0 w)\leq\sup_{\tau>0}\mathcal{I}_{V(x_0)}(\tau w)&=\mathcal{I}_{V(x_0)}(w)\\
&=c_{V(x_0)}
\end{align*}
which contradicts with \eqref{equ6-4}. Thus $V(x_0)=V_0$.
\end{proof}

From $V(x_0)=V_0$, we see that $c_{V(x_0)}=c_{V_0}$. It follows from \eqref{equ6-4} that $\|v_n\|\rightarrow\|v\|$. Using Brezis-Lieb Lemma, we conclude that $v_n\rightarrow v$ in $H^s(\mathbb{R}^3)$.

\begin{lemma}\label{lem6-4}
Let $v_n(z)=u_n(z+y_n)$ satisfies \eqref{equ6-2}. Then
\begin{equation*}
\lim_{|z|\rightarrow\infty}v_n(z)=0\quad \text{uniformly for}\,\, n\in\mathbb{N}.
\end{equation*}
\end{lemma}

\begin{proof}
By Proposition \ref{pro5-1}, we see that there exists $C>0$ independent of $n$ such that $\|v_n\|_{L^{\infty}}\leq C$. Now, we rewrite the reduced form of problem \eqref{equ6-1} as follows
\begin{equation*}
(-\Delta)^sv_n+v_n=h_n(z)\quad z\in\mathbb{R}^3,
\end{equation*}
where $h_n(z):=v_n-V(\varepsilon_n z+\varepsilon_ny_n)v_n-\phi_{v_n}^tv_n+g(\varepsilon_n z+\varepsilon_ny_n,v_n)$.
Clearly, $h_n\in L^{\infty}(\mathbb{R}^3)$ and is uniformly bounded. By interpolation inequality and $v_n\rightarrow v$ in $H^s(\mathbb{R}^3)$, for $n\rightarrow\infty$, we have that $h_n\rightarrow h$ in $L^r(\mathbb{R}^3)$ for $2\leq r<+\infty$, where $h(z)=v(z)-V(x_0)v(z)-\phi_{v}^t(z)v(z)+g(x_0,v(z))$. Using some results found in \cite{FQT}, we see that
\begin{equation*}
v_n(z)=\int_{\mathbb{R}^3}\mathcal{K}(z-y)h_n(y)\,{\rm d}y
\end{equation*}
where $\mathcal{K}$ is a Bessel potential, which possesses the following properties:\\
$(\mathcal{K}_1)$ $\mathcal{K}$ is positive, radially symmetric and smooth in $\mathbb{R}^3\backslash\{0\}$;\\
$(\mathcal{K}_2)$ there exists a constant $C>0$ such that $\mathcal{K}(x)\leq \frac{C}{|x|^{3+2s}}$ for all $x\in\mathbb{R}^3\backslash\{0\}$;\\
$(\mathcal{K}_3)$ $\mathcal{K}\in L^\tau(\mathbb{R}^3)$ for $\tau\in[1,\frac{3}{3-2s})$.

We define two sets $A_{\delta}=\{y\in\mathbb{R}^3\,\,|\,\, |z-y|\geq\frac{1}{\delta}\}$ and $B_{\delta}=\{y\in\mathbb{R}^3\,\,|\,\, |z-y|<\frac{1}{\delta}\}$.
Hence,
\begin{equation*}
0\leq v_n(z)\leq \int_{\mathbb{R}^3}\mathcal{K}(z-y)|h_n(y)|\,{\rm d}y=\int_{A_{\delta}}\mathcal{K}(z-y)|h_n(y)|\,{\rm d}y+\int_{B_{\delta}}\mathcal{K}(z-y)|h_n(y)|\,{\rm d}y.
\end{equation*}
From the definition of $A_{\delta}$ and $(\mathcal{K}_2)$, we have that for all $n\in\mathbb{N}$,
\begin{equation*}
\int_{A_{\delta}}\mathcal{K}(z-y)|h_n(y)|\,{\rm d}y\leq C\delta^s\|h_n\|_{\infty}\int_{A_{\delta}}\frac{1}{|z-y|^{3+s}}\,{\rm d}y\leq C\delta^s\int_{A_{\delta}}\frac{1}{|z-y|^{3+s}}\,{\rm d}y:=C\delta^{2s}.
\end{equation*}

On the other hand, by H\"{o}lder's inequality and $(\mathcal{K}_3)$, we deduce that
\begin{align*}
&\int_{B_{\delta}}\mathcal{K}(z-y)|h_n(y)|\,{\rm d}y\leq\int_{B_{\delta}}\mathcal{K}(z-y)|h_n-h|\,{\rm d}y+\int_{B_{\delta}}\mathcal{K}(z-y)|h|\,{\rm d}y\\
&\leq\Big(\int_{B_{\delta}}\mathcal{K}^{\frac{6}{3+2s}}\,{\rm d}y\Big)^{\frac{3+2s}{6}}\Big(\int_{B_{\delta}}|h_n-h|^{\frac{6}{3-2s}}\,{\rm d}y\Big)^{\frac{3-2s}{6}}+
\Big(\int_{B_{\delta}}\mathcal{K}^2\,{\rm d}y\Big)^{\frac{1}{2}}\Big(\int_{B_{\delta}}|h|^2\,{\rm d}y\Big)^{\frac{1}{2}}\\
&\leq\Big(\int_{\mathbb{R}^3}\mathcal{K}^{\frac{6}{3+2s}}\,{\rm d}y\Big)^{\frac{3+2s}{6}}\Big(\int_{\mathbb{R}^3}|h_n-h|^{\frac{6}{3-2s}}\,{\rm d}y\Big)^{\frac{3-2s}{6}}+
\Big(\int_{\mathbb{R}^3}\mathcal{K}^2\,{\rm d}y\Big)^{\frac{1}{2}}\Big(\int_{B_{\delta}}|h|^2\,{\rm d}y\Big)^{\frac{1}{2}}
\end{align*}
where we have used the fact that $s>\frac{3}{4}$ so that $\frac{6}{3+2s}<\frac{3}{3-2s}$ and $2<\frac{3}{3-2s}$.

Since $\Big(\int_{B_{\delta}}|h|^2\,{\rm d}y\Big)^{\frac{1}{2}}\rightarrow0$ as $|z|\rightarrow+\infty$, thus, we deduce that there exist $n_0\in\mathbb{N}$ and $R_0>0$ independence of $\delta>0$ such that
\begin{equation*}
\int_{B_{\delta}}\mathcal{K}(z-y)|h_n(y)|\,{\rm d}y\leq \delta,\quad \forall n\geq n_0\quad \text{and}\quad |z|\geq R_0.
\end{equation*}
Hence,
\begin{equation*}
\int_{\mathbb{R}^3}\mathcal{K}(z-y)|h_n(y)|\,{\rm d}y\leq C\delta^{2s}+\delta,\quad \forall n\geq n_0\quad \text{and}\quad |z|\geq R_0.
\end{equation*}
For each $n\in\{1,2,\cdots,n_0-1\}$, there exists $R_n>0$ such that $\Big(\int_{B_{\delta}}|h_n|^2\,{\rm d}y\Big)^{\frac{1}{2}}< \delta$ as $|z|\geq R_n$. Thus, for $|z|\geq R_n$, we have that
\begin{align*}
\int_{\mathbb{R}^3}\mathcal{K}(z-y)|h_n(y)|\,{\rm d}y&\leq C\delta^{2s}+\int_{B_{\delta}}\mathcal{K}(z-y)|h_n(y)|\,{\rm d}y\\
&\leq C\delta^{2s}+\|\mathcal{K}\|_2\Big(\int_{B_{\delta}}|h_n|^2\,{\rm d}y\Big)^{\frac{1}{2}}\leq C(\delta^{2s}+\delta)
\end{align*}
for each $n\in\{1,2,\cdots,n_0-1\}$.
Therefore, taking $R=\max\{R_0,R_1,\cdots, R_{n_0-1}\}$, we infer that for any $n\in\mathbb{N}$, there holds
\begin{equation*}
0\leq v_n(z)\leq\int_{\mathbb{R}^3}\mathcal{K}(z-y)|h_n(y)|\,{\rm d}y\leq C\delta^{2s}+\delta,\quad \text{for all}\quad |z|\geq R
\end{equation*}
implies that $\lim\limits_{|z|\rightarrow\infty}v_n(z)=0$ uniformly in $n\in\mathbb{N}$.
\end{proof}

\begin{lemma}\label{lem6-5}
There is $n_0>0$ such that $u_n(z)=v_n(z-y_n)<a$, for all $n\geq n_0$ and all $z\in\mathbb{R}^3\backslash(\Lambda/\varepsilon_n)$. Hence, $v_n$ is a solution of problem \eqref{equ3-1} for $n\geq n_0$.
\end{lemma}
\begin{proof}
By Lemma \ref{lem6-3}, we see that $\varepsilon_ny_n\rightarrow x_0$ and $x_0\in\Lambda$. Thus, there exists $R'>0$ such that for some subsequence, still denoted by itself, $B_{R'}(\varepsilon_ny_n)\subset\Lambda\subset \Lambda'$ for all $n\in\mathbb{N}$. Hence, $B_{R'/\varepsilon_n}(y_n)\subset\Lambda/\varepsilon_n$, $\forall n\in\mathbb{N}$. Moreover, by Lemma \ref{lem6-4}, there is $R_1>0$ such that $v_n(z)<a$ for $|z|\geq R_1$ and $\forall n\in\mathbb{N}$.  Thus,
\begin{equation*}
u_n(z)=v_n(z-y_n)<a,\quad \text{for all}\,\,z\in \mathbb{R}^3\backslash B_{R_1}(y_n)\,\,\text{and}\,\,\forall n\in\mathbb{N}
\end{equation*}
Hence, there exists $n_0\in\mathbb{N}$ such that
\begin{equation*}
\mathbb{R}^3\backslash(\Lambda'/\varepsilon_n)\subset\mathbb{R}^3\backslash(\Lambda/\varepsilon_n)\subset\mathbb{R}^3\backslash B_{R'/\varepsilon_n}(y_n)\subset\mathbb{R}^3\backslash B_{R_1}(y_n), \quad \forall n\geq n_0
\end{equation*}
and then
\begin{equation*}
u_n(z)<a\quad \forall z\in\mathbb{R}^3\backslash(\Lambda/\varepsilon_n)\,\, \text{and} \,\, \forall\, n\geq n_0.
\end{equation*}
\end{proof}

{\bf Proof of Theorem \ref{thm1-1}.}
By Proposition \ref{pro3-1}, we see that problem \eqref{equ3-1} has a nonnegative solution $v_{\varepsilon}$ for all $\varepsilon>0$. From Lemma \ref{lem6-5}, there exists $\varepsilon_0>0$ such that
\begin{equation}\label{equ6-5}
v_{\varepsilon }(z)<a\quad \forall z\in\mathbb{R}^3\backslash(\Lambda/\varepsilon)\,\,\text{and}\,\, \varepsilon\in(0,\varepsilon_0)
\end{equation}
which implies that $g(\varepsilon z,v_{\varepsilon})=f(v_{\varepsilon})+v_{\varepsilon}^{2_s^{\ast}-1}$. Thus, $v_{\varepsilon}$ is a solution of problem
\begin{equation}\label{equ6-6}
(-\Delta)^sv+V(\varepsilon z) v+ \phi_{v}^t v=f(v)+v^{2_s^{\ast}-1}\quad z\in\mathbb{R}^3.
\end{equation}
for all $\varepsilon\in(0,\varepsilon_0)$. Let $u_{\varepsilon}(x)=v_{\varepsilon}(x/\varepsilon)$ for every $\varepsilon\in(0,\varepsilon_0)$, it follows that $u_{\varepsilon}$ must be a solution to original problem \eqref{main} for $\varepsilon\in(0,\varepsilon_0)$.

If $z_{\varepsilon}$ denotes a global maximum point of $v_{\varepsilon}$, then
\begin{equation}\label{equ6-7}
v_{\varepsilon}(z_{\varepsilon})\geq a\quad \forall\, \varepsilon\in(0,\varepsilon_0).
\end{equation}
Suppose that $v_{\varepsilon}(z_{\varepsilon})<a$, taking $v_{\varepsilon}$ as a text function for \eqref{equ6-6}, we get
\begin{align*}
\alpha_0\int_{\mathbb{R}^3}v_{\varepsilon}^2\,{\rm d}z&\leq\int_{\mathbb{R}^3}V(\varepsilon z)v_{\varepsilon}^2\,{\rm d}z\leq\int_{\mathbb{R}^3}(f(v_{\varepsilon})v_{\varepsilon}+v_{\varepsilon}^{2_s^{\ast}})\,{\rm d}z\\
&=\int_{\mathbb{R}^3}v_{\varepsilon}^2\Big(\frac{f(v_{\varepsilon})}{v_{\varepsilon}^3}v_{\varepsilon}^2+v_{\varepsilon}^{2_s^{\ast}-2}\Big)\,{\rm d}z\\
&\leq\int_{\mathbb{R}^3}v_{\varepsilon}^2(\frac{f(a)}{a}+a^{2_s^{\ast}-2})\,{\rm d}z=\frac{\alpha_0}{k}\int_{\mathbb{R}^3}v_{\varepsilon}^2\,{\rm d}z.
\end{align*}
Hence we get a contradiction owing to the choosing $k>2$. In view of Lemma \ref{lem6-4}, we see that $\{z_{\varepsilon}\}$ is bounded for $\varepsilon\in(0,\varepsilon_0)$.

In what follows, setting $x_{\varepsilon}=\varepsilon z_{\varepsilon}+\varepsilon y_{\varepsilon}$, where $\{y_{\varepsilon}\}$ is given in Lemma \ref{lem6-1}. Since $u_{\varepsilon}(x)=v_{\varepsilon}(\frac{x}{\varepsilon}-y_{\varepsilon})$, then $x_{\varepsilon}$ is a global maximum point of $u_{\varepsilon}$ and $u_{\varepsilon}(x_{\varepsilon})\geq a$ for all $\varepsilon\in(0,\varepsilon_0)$.

Now, we claim that $\lim\limits_{\varepsilon\rightarrow0^{+}}V(x_{\varepsilon})=V_0$. Indeed, if the above limit does not hold, there is $\varepsilon_n\rightarrow0^{+}$ and $\gamma_0>0$ such that
\begin{equation}\label{equ6-8}
V(x_{\varepsilon_n})\geq V_0+\gamma_0\quad \forall n\in\mathbb{N}.
\end{equation}

By Lemma \ref{lem6-4}, we know that $\lim\limits_{|z|\rightarrow\infty}v_{\varepsilon_n}(z)=0$ uniformly in $n\in\mathbb{N}$. From \eqref{equ6-7}, thus $\{z_{\varepsilon_n}\}$ is a bounded sequence. Up to a subsequence, using Lemma \ref{lem6-3}, we know that there is $x_0\in\Lambda$ such that $V(x_0)=V_0$ and $\varepsilon_ny_{\varepsilon_n}\rightarrow x_0$. Hence, $x_{\varepsilon_n}=\varepsilon_n z_{\varepsilon_n}+\varepsilon_n y_{\varepsilon_n}\rightarrow x_0$ which implies that $V(x_{\varepsilon_n})\rightarrow V(x_0)=V_0$ contradicting with \eqref{equ6-8}.

To complete the proof, we only need to prove the decay properties of $u_{\varepsilon}$. Similar argument to the proof of Lemma 5.6 in \cite{Teng2}, we can obtain that
\begin{equation*}
0<v_{\varepsilon}(z)\leq \frac{C}{1+|z|^{3+2s}}.
\end{equation*}
Thus, by the boundedness of $\{z_{\varepsilon}\}$, i.e., there exists $C_0>0$ such that $|z_{\varepsilon}|\leq C_0$, we have
\begin{align*}
u_{\varepsilon}(x)=v_{\varepsilon}(\frac{x}{\varepsilon}-y_{\varepsilon})\leq \frac{C}{1+|\frac{x-x_{\varepsilon}+\varepsilon z_{\varepsilon}}{\varepsilon}|^{3+2s}}&\leq\frac{ C\varepsilon^{3+2s}}{\varepsilon^{3+2s}(1-C_0^{3+2s})+|x-x_{\varepsilon}|^{3+2s}}\\
&:=\frac{ C\varepsilon^{3+2s}}{\varepsilon^{3+2s}C_1+|x-x_{\varepsilon}|^{3+2s}}.
\end{align*}

\section{Multiplicity of solutions to \eqref{main}}
In this section, we will use the following two abstract Propositions to get the multiplicity of solutions.
\begin{proposition}\label{pro7-1}(\cite{ZG})
Let $I$ be a $C^1$-functional defined on a $C^1$-Finsler manifold $M$. If $I$ is bounded from below and satisfies the $(PS)$ condition, then $I$ possesses at least ${\rm cat}_{M}(M)$ distinct critical point.
\end{proposition}

Let us consider $\delta>0$ such that $\mathcal{M}_{\delta}\subset\Lambda$ and a smooth cut-off function with $0\leq\eta\leq1$, $\eta=1$ on $B_{\delta/2}(0)$, $\eta=0$ on $\mathbb{R}^3\backslash B_{\delta}(0)$, $|\nabla\eta|\leq C$. For any $y\in \mathcal{M}$, set
\begin{equation*}
\psi_{\varepsilon,y}(z)=\eta(\varepsilon z-y)w(\frac{\varepsilon z-y}{\varepsilon}),
\end{equation*}
where $w\in H^s(\mathbb{R}^3)$ is a solution of \eqref{equ4-1} with $\mu=V_0$ such that $\mathcal{I}_{V_0}'(w)=0$ and $\mathcal{I}_{V_0}=c_{V_0}$. Thus, there exists $t_{\varepsilon}>0$ such that $\max\limits_{t>0}J_{\varepsilon}(t\psi_{\varepsilon,y})=J_{\varepsilon}(t_{\varepsilon}\psi_{\varepsilon,y})$. We define $\Phi_{\varepsilon}:\mathcal{M}\rightarrow\mathcal{N}_{\varepsilon}$ by
\begin{equation*}
\Phi_{\varepsilon}(y)=t_{\varepsilon}\psi_{\varepsilon,y}
\end{equation*}

For the $\delta>0$ given by above, we choose $\rho=\rho(\delta)>0$ such that $\mathcal{M}_{\delta}\subset B_{\rho}(0)$. Define $\Upsilon:\mathbb{R}^3\rightarrow\mathbb{R}$ as $\Upsilon(z)=z$ for $|z|\leq\rho$ and $\Upsilon(z)=\frac{\rho z}{|z|}$ for $|z|\geq\rho$, and consider the map $\beta_{\varepsilon}:\mathcal{N}_{\varepsilon}\rightarrow\mathbb{R}^3$ given by
\begin{equation*}
\beta_{\varepsilon}(u)=\frac{\int_{\mathbb{R}^3}\Upsilon(\varepsilon z)u^2\,{\rm d}z}{\int_{\mathbb{R}^3}u^2\,{\rm d}z}
\end{equation*}

Define
\begin{equation*}
\widetilde{\mathcal{N}_{\varepsilon}}=\{u\in\mathcal{N}_{\varepsilon}\,\,\Big|\,\,J_{\varepsilon}(u)\leq c_{V_0}+h(\varepsilon)\},
\end{equation*}
where $h(\varepsilon)=\sup\limits_{y\in\mathcal{M}}|J_{\varepsilon}(\Phi_{\varepsilon}(y))-c_{V_0}|$.

\begin{proposition}\label{pro7-2}(\cite{BC0}, Lemma 4.3)
Let $\Phi_{\varepsilon}:\mathcal{M}\rightarrow \widetilde{\mathcal{N}_{\varepsilon}}$, $\beta_{\varepsilon}:\widetilde{\mathcal{N}_{\varepsilon}}\rightarrow \mathcal{M}_{\delta}$ be two continuous maps defined as above. If $\beta_{\varepsilon}\circ\Phi_{\varepsilon}$ is homotopically equivalent to the embedding ${\rm id}:\mathcal{M}\rightarrow \mathcal{M}_{\delta}$. Then ${\rm cat}_{\widetilde{\mathcal{N}_{\varepsilon}}}(\widetilde{\mathcal{N}_{\varepsilon}})\geq{\rm cat}_{\mathcal{M}_{\delta}}(\mathcal{M})$.
\end{proposition}

Therefore, in Proposition \ref{pro7-1}, we choose the Finsler manifold $M$ as $\widetilde{\mathcal{N}_{\varepsilon}}$. It is standard to show the following result.
\begin{proposition}
For any $\delta>0$, there exists $\varepsilon_{\delta}>0$ such that for any $\varepsilon\in(0,\varepsilon_{\delta})$, the system \eqref{equ3-1} has at least ${\rm cat}_{\mathcal{M}_{\delta}}(\mathcal{M})$ solutions, where $\mathcal{M}$ and $\mathcal{M}_{\delta}$ defined in Introduction.
\end{proposition}

The remain is to verify the $(PS)$ condition and the homopotically equivalent of $\beta_{\varepsilon}\circ\Phi_{\varepsilon}$ with the embedding ${\rm id}:\mathcal{M}\rightarrow\mathcal{M}_{\delta}$. The proof is standard, we are only to verify the $(PS)$ condition. The other detailed proof can be consulted in the papers \cite{ASA,AY,HLP} and the references therein.

\begin{lemma}\label{lem7-1}
Let $\{u_n\}\subset H_{\varepsilon}$ be a $(PS)_c$ sequence for $c\in(0,\frac{s}{3}\mathcal{S}_s^{\frac{3}{2s}})$. Then for each $\delta>0$, there exists $R>0$ such that
\begin{equation}\label{equ7-1}
\limsup_{n\rightarrow\infty}\int_{\mathbb{R}^3\backslash B_R(0)}(|(-\Delta)^{\frac{s}{2}}u_n|^2+V(\varepsilon z)u_n^2)\,{\rm d}z<\delta.
\end{equation}
\end{lemma}
\begin{proof}
From Lemma \ref{lem3-3}, we know that $\{u_n\}$ is bounded in $H_{\varepsilon}$ and up to a subsequence, we may assume that there exists $u\in H_{\varepsilon}$ such that $u_n\rightharpoonup u$ in $H_{\varepsilon}$, $u_n\rightarrow u$ in $L_{loc}^r(\mathbb{R}^3)$ for $1\leq r<2_s^{\ast}$ and $u_n\rightarrow u$ a.e. $\mathbb{R}^3$. First we may assume that $R$ is chosen such that $\Lambda'/\varepsilon\subset B_{R/2}(0)$. Let $\eta_R$ be a smooth cut-off function so that $\eta_R=0$ on $B_{R/2}(0)$, $\eta_R=1$ on $\mathbb{R}^3\backslash B_R(0)$, $0\leq\eta_R\leq1$ and $|\nabla\eta_R|\leq\frac{C}{R}$. Since $\{u_n\}$ is a bounded $(PS)_c$ sequence, we have
\begin{align*}
&\int_{\mathbb{R}^3}(-\Delta)^{\frac{s}{2}}u_n(-\Delta)^{\frac{s}{2}}(\eta_Ru_n)\,{\rm d}z+\int_{\mathbb{R}^3}V(\varepsilon z)u_n^2\eta_R\,{\rm d}z+\int_{\mathbb{R}^3}\phi_{u_n}^tu_n^2\eta_R\,{\rm d}z\\
&=\int_{\mathbb{R}^3}g(\varepsilon z,u_n)u_n\eta_R\,{\rm d}z+o(1)\leq\frac{1}{k}\int_{\mathbb{R}^3}V(\varepsilon z)u_n^2\eta_R\,{\rm d}z+o(1)
\end{align*}
Similar arguments to \eqref{equ2-5} and \eqref{equ2-6}, we deduce that
\begin{align}\label{equ7-1-0}
&\int_{\mathbb{R}^3\backslash B_R(0)}(|(-\Delta)^{\frac{s}{2}}u_n|^2+(1-\frac{1}{k})V(\varepsilon z)u_n^2)\,{\rm d}z\nonumber\\
&\leq\int_{\mathbb{R}^3}(|(-\Delta)^{\frac{s}{2}}u_n|^2+(1-\frac{1}{k})V(\varepsilon z)u_n^2)\eta_R\,{\rm d}z\leq\frac{C}{R^{\frac{s}{2}}}\|u_n\|_{H_{\varepsilon}}^2+o(1)
\end{align}
which implies that \eqref{equ7-1} holds.
\end{proof}

By the well known argument, we see that the Nehari manifold $\mathcal{N}_{\varepsilon}$ is a $C^1$-manifold.
\begin{lemma}\label{lem7-2}
The functional $J_{\varepsilon}$ restricted to $\mathcal{N}_{\varepsilon}$ satisfies $(PS)_c$ condition for each $c\in(0,\frac{s}{3}\mathcal{S}_s^{\frac{3}{2s}})$.
\end{lemma}
\begin{proof}
Let $\{u_n\}\subset\mathcal{N}_{\varepsilon}$ be such that
\begin{equation*}
J_{\varepsilon}(u_n)\rightarrow c\quad \text{and}\quad \|J_{\varepsilon}'(u_n)\|_{\ast}\rightarrow0\quad \text{as}\,\, n\rightarrow\infty,
\end{equation*}
where $\|J_{\varepsilon}'(u)$ denotes the norm of the derivative of $J_{\varepsilon}$ restricted to $\mathcal{N}_{\varepsilon}$ at the point $u\in\mathcal{N}_{\varepsilon}$. Similar arguments to the proof of Lemma \ref{lem3-3}, we obtain that $\{u_n\}$ is bounded in $H_{\varepsilon}$. Thus, up to a subsequence, we may assume that there is $u\in H_{\varepsilon}$ such that
\begin{equation}\label{equ7-2}
\left\{
  \begin{array}{ll}
    u_n\rightharpoonup u& \hbox{in $H_{\varepsilon}$,} \\
    u_n\rightarrow u & \hbox{in $L_{loc}^r(\mathbb{R}^3)$ for $1\leq r<2_s^{\ast}$,} \\
    u_n\rightarrow u & \hbox{a.e. $\mathbb{R}^3$.}
  \end{array}
\right.
\end{equation}
By standard computation, we can assume that there exists $\{\lambda_n\}\subset\mathbb{R}$ such that
\begin{equation*}
J_{\varepsilon}'(u_n)=\lambda_n\mathcal{G}_{\varepsilon}'(u_n)+o(1),
\end{equation*}
where $\mathcal{G}_{\varepsilon}(u)=\langle J_{\varepsilon}'(u),u\rangle$. Moreover, since $u_n\in\mathcal{N}_{\varepsilon}$, we know that
\begin{equation}\label{equ7-3}
0=\langle J_{\varepsilon}'(u_n),u_n\rangle=\lambda_n\langle\mathcal{G}_{\varepsilon}'(u_n),u_n\rangle+o(1)\|u_n\|_{H_{\varepsilon}}.
\end{equation}
Next we will show that $\lambda_n\rightarrow0$ as $n\rightarrow\infty$. Indeed, by the fact that $(\frac{f(\tau)}{\tau^3})'>0$, $(\frac{\tilde{f}(\tau)}{\tau})'>0$ for all $\tau\geq0$, we deduce that
\begin{align*}
&\langle\mathcal{G}_{\varepsilon}'(u_n),u_n\rangle\\
&=2\int_{\mathbb{R}^3}(|(-\Delta)^{\frac{s}{2}}u_n|^2+V(\varepsilon z)u_n^2)\,{\rm d}z+4\int_{\mathbb{R}^3}\phi_{u_n}^tu_n^2\,{\rm d}z-\int_{\mathbb{R}^3}(g'(\varepsilon z,u_n)u_n^2\\
&+g(\varepsilon z,u_n)u_n)\,{\rm d}z\\
&=-2\int_{\mathbb{R}^3}(|(-\Delta)^{\frac{s}{2}}u_n|^2+V(\varepsilon z)u_n^2)\,{\rm d}z+\int_{\mathbb{R}^3}(3g(\varepsilon z,u_n)u_n-g'(\varepsilon z,u_n)u_n^2)\,{\rm d}z\\
&\leq-2\int_{\mathbb{R}^3}(|(-\Delta)^{\frac{s}{2}}u_n|^2+V(\varepsilon z)u_n^2)\,{\rm d}z+2\int_{\mathbb{R}^3}(1-\chi(\varepsilon z)\tilde{f}(u_n)u_n\,{\rm d}z\\
&\leq-2\int_{\mathbb{R}^3}(|(-\Delta)^{\frac{s}{2}}u_n|^2+V(\varepsilon z)u_n^2)\,{\rm d}z+\frac{2}{k}\int_{\mathbb{R}^3}V(\varepsilon z)u_n^2\,{\rm d}z\\
&\leq-2(1-\frac{1}{k})\|u_n\|_{H_{\varepsilon}}^2.
\end{align*}
Thus, we may assume that $\langle\mathcal{G}_{\varepsilon}'(u_n),u_n\rangle\rightarrow l<0$. It follows from \eqref{equ7-3} that $\lambda_n\rightarrow0$ as $n\rightarrow\infty$ and then we see that $J_{\varepsilon}'(u_n)\rightarrow0$ as $n\rightarrow\infty$ in the dual space of $H_{\varepsilon}$. Hence, $\{u_n\}$ is a $(PS)_c$ sequence for $J_{\varepsilon}$.

We claim that
\begin{equation}\label{equ7-4}
\int_{\mathbb{R}^3}g(\varepsilon z,u_n)u_n\,{\rm d}z\rightarrow\int_{\mathbb{R}^3}g(\varepsilon z,u)u\,{\rm d}z\quad \text{as}\,\, n\rightarrow\infty.
\end{equation}
In fact, by \eqref{equ7-1-0}, we can also obtain that
\begin{equation}\label{equ7-5}
\int_{\mathbb{R}^3}u_n^2\eta_R^2\,{\rm d}z\leq\frac{C}{R^{\frac{s}{2}}}+o(1)\quad \int_{\mathbb{R}^3}|(-\Delta)^{\frac{s}{2}}u_n|^2\eta_R^2\,{\rm d}z\leq\frac{C}{R^{\frac{s}{2}}}+o(1).
\end{equation}
By interpolation inequality, we have that for any $\frac{1}{r}=\frac{\theta}{2}+\frac{1-\theta}{2_s^{\ast}}$ with $0<\theta\leq1$,
\begin{equation*}
\Big(\int_{\mathbb{R}^3}|u_n\eta_R|^r\,{\rm d}z\Big)^{\frac{1}{r}}\leq\Big(\int_{\mathbb{R}^3}u_n^2\eta_R^2\,{\rm d}z\Big)^{\frac{\theta}{2}}\Big(\int_{\mathbb{R}^3}|u_n\eta_R|^{2_s^{\ast}}\,{\rm d}z\Big)^{\frac{1-\theta}{2_s^{\ast}}}\leq \frac{C}{R^{\frac{\theta s}{4}}}+o(1)
\end{equation*}
which yields
\begin{equation*}
\int_{\mathbb{R}^3\backslash B_R(0)}|u_n|^r\,{\rm d}z\leq \frac{C}{R^{\frac{r\theta s}{4}}}+o(1)\quad \text{for all}\,\, 2\leq r<2_s^{\ast}.
\end{equation*}
By Sobolev inequality, \eqref{equ7-5} and using similar arguments to Lemma \ref{lem3-5}, we have that
\begin{align*}
\Big(\int_{\mathbb{R}^3}|u_n\eta_R|^{2_s^{\ast}}\,{\rm d}z\Big)^{\frac{2}{2_s^{\ast}}}\leq\mathcal{S}_s^{-1}\int_{\mathbb{R}^3}|(-\Delta)^{\frac{s}{2}}u_n|^2\eta_R^2\,{\rm d}z+\frac{C}{R^{\frac{s}{2}}}\leq\frac{C}{R^{\frac{s}{2}}}+o(1).
\end{align*}
Therefore, we get
\begin{equation}\label{equ7-6}
\int_{\mathbb{R}^3\backslash B_R(0)}|u_n|^r\,{\rm d}z\leq\frac{C}{R^{\frac{s}{2}}}+o(1)\quad \text{for all}\,\, 2\leq r\leq 2_s^{\ast},
\end{equation}
which together with \eqref{equ7-2} implies that
\begin{equation}\label{equ7-7}
u_n\rightarrow u\quad \text{in}\,\, L^r(\mathbb{R}^3)\,\, \text{for all}\,\, 2\leq r<2_s^{\ast}.
\end{equation}
By the estimate \eqref{equ7-6}, it is easy to check that
\begin{equation*}
\int_{\mathbb{R}^3\backslash B_R(0)}g(\varepsilon z,u_n)u_n\,{\rm d}z\leq\frac{C}{R^{\frac{s}{2}}}+o(1).
\end{equation*}
On the other hand, using \eqref{equ7-2} and Lebesgue dominated convergence theorem, it is easy to show that
\begin{equation*}
\int_{B_R(0)}\chi(\varepsilon z)f(u_n)u_n\,{\rm d}z\rightarrow\int_{B_R(0)}\chi(\varepsilon z)f(u)u\,{\rm d}z
\end{equation*}
and
\begin{equation*}
\int_{B_R(0)}(1-\chi(\varepsilon z))\tilde{f}(u_n)u_n\,{\rm d}z\rightarrow\int_{B_R(0)}(1-\chi(\varepsilon z))\tilde{f}(u)u\,{\rm d}z.
\end{equation*}
In order to prove \eqref{equ7-4}, it is only need to show that
\begin{equation}\label{equ7-8}
\int_{B_R(0)}\chi(\varepsilon z)(u_n^{+})^{2_s^{\ast}}\,{\rm d}z\rightarrow\int_{B_R(0)}\chi(\varepsilon z)(u^{+})^{2_s^{\ast}}\,{\rm d}z\quad \text{as}\,\, n\rightarrow\infty.
\end{equation}
Now we apply Lemma \ref{lem2-2} to establish \eqref{equ7-8}. Since $\{u_n\}$ is bounded in $H^s(\mathbb{R}^3)$, by Phrokorov¡¯s theorem (Theorem 8.6.2 in \cite{B}), there exist $\mu,\nu\in\mathcal{M}(\mathbb{R}^3)$ such that
\begin{equation*}
|(-\Delta)^{\frac{s}{2}}u_n^{+}|^2\rightharpoonup\mu\,\,\text{and}\,\, (u_n^{+})^{2_s^{\ast}}\rightharpoonup\nu\,\, \text{weakly-}\ast\,\, \text{in}\,\, \mathcal{M}(\mathbb{R}^3)\,\,\text{as}\,\, n\rightarrow\infty.
\end{equation*}
By \eqref{equ7-7} and using Lemma \ref{lem2-2}, there exist an at most countable index set $J$, sequence $\{x_j\}_{j\in J}\subset\mathbb{R}^3$ and $\{\mu_j\}$, $\{\nu_j\}\subset(0,\infty)$ such that
\begin{equation}\label{equ7-9}
\nu=(u^{+})^{2_s^{\ast}}+\sum_{j\in J}\nu_j\delta_{x_j},
\end{equation}
\begin{equation}\label{equ7-10}
\limsup_{n\rightarrow\infty}\int_{\mathbb{R}^3}|(-\Delta)^{\frac{s}{2}}u_n^{+}|^2\,{\rm d}z=\int_{\mathbb{R}^3}\,{\rm d}\mu+\mu_{\infty},\quad \limsup_{n\rightarrow\infty}\int_{\mathbb{R}^3}(u_n^{+})^{2_s^{\ast}}\,{\rm d}z=\int_{\mathbb{R}^3}\,{\rm d}\nu+\nu_{\infty}
\end{equation}
and
\begin{equation}\label{equ7-11}
\nu_j\leq(\mathcal{S}_s^{-1}\mu(\{x_j\}))^{\frac{2_s^{\ast}}{2}}\,\, \text{for any}\,\, j\in J\,\, \text{and}\,\,\nu_{\infty}\leq(\mathcal{S}_s^{-1}\mu_{\infty})^{\frac{2_s^{\ast}}{2}}.
\end{equation}
It is suffices to show that $\{x_j\}_{j\in J}\cap\{z|\chi(\varepsilon z)>0\}=\emptyset$. Suppose by contradiction that $\chi(\varepsilon x_j)>0$ for some $j\in J$. Define the function $\psi_{\rho}(z)=\psi(\frac{z-x_j}{\rho})$ for $\rho>0$, where $\psi$ is a smooth cut-off function such that $\psi=1$ on $B_1(0)$, $\psi=0$ on $\mathbb{R}^3\backslash B_2(0)$, $0\leq\psi\leq1$ and $|\nabla\psi|\leq C$. Suppose that $\rho$ is chosen in such a way that the support of $\psi_{\rho}$ is contained in $\{z|\chi(\varepsilon z)>0\}$.
Since
\begin{equation*}
\langle J_{\varepsilon}'(u_n),\psi_{\rho}u_n^{+}\rangle\rightarrow0\quad \text{as}\,\, n\rightarrow\infty,
\end{equation*}
i.e.,
\begin{align}\label{equ7-12}
\int_{\mathbb{R}^3}(-\Delta)^{\frac{s}{2}}u_n(-\Delta)^{\frac{s}{2}}(u_n^{+}\psi_{\rho})\,{\rm d}z&+\int_{\mathbb{R}^3}V(\varepsilon z)(u_n^{+})^2\psi_{\rho}\,{\rm d}z+\int_{\mathbb{R}^3}\phi_{u_n}^t(u_n^{+})^2\psi_{\rho}\,{\rm d}z\nonumber\\
&=\int_{\mathbb{R}^3}g(\varepsilon z,u_n)u_n^{+}\psi_{\rho}\,{\rm d}z+o(1).
\end{align}
By \eqref{equ7-2}, we can deduce that
\begin{align*}
\lim_{\rho\rightarrow0^{+}}\limsup_{n\rightarrow\infty}\int_{B_{2\rho}(x_j)}V(\varepsilon z)(u_n^{+})^2\psi_{\rho}\,{\rm d}z=\lim_{\rho\rightarrow0^{+}}\int_{B_{2\rho}(x_j)}V(\varepsilon z)(u^{+})^2\psi_{\rho}\,{\rm d}z=0,
\end{align*}
\begin{align*}
\lim_{\rho\rightarrow0^{+}}\limsup_{n\rightarrow\infty}\int_{B_{2\rho}(x_j)}\phi_{u_n}^t(u_n^{+})^2\psi_{\rho}\,{\rm d}z=\lim_{\rho\rightarrow0^{+}}\int_{B_{2\rho}(x_j)}\phi_u^t(u^{+})^2\psi_{\rho}\,{\rm d}z=0,
\end{align*}
and similarly,
\begin{equation*}
\lim_{\rho\rightarrow0^{+}}\limsup_{n\rightarrow\infty}\int_{B_{2\rho}(x_j)}\chi(\varepsilon z)f(u_n)u_n^{+}\psi_{\rho}\,{\rm d}z=0,
\end{equation*}
\begin{equation*}
\lim_{\rho\rightarrow0^{+}}\limsup_{n\rightarrow\infty}\int_{B_{2\rho}(x_j)}(1-\chi(\varepsilon z))\tilde{f}(u_n)u_n^{+}\psi_{\rho}\,{\rm d}z=0
\end{equation*}
which leads to
\begin{align*}
&\lim_{\rho\rightarrow0^{+}}\limsup_{n\rightarrow\infty}\int_{B_{2\rho}(x_j)}g(\varepsilon z,u_n)u_n^{+}\psi_{\rho}\,{\rm d}z=\lim_{\rho\rightarrow0^{+}}\limsup_{n\rightarrow\infty}\int_{B_{2\rho}(x_j)}\chi(\varepsilon z)(u_n^{+})^{2_s^{\ast}}\psi_{\rho}\,{\rm d}z\\
&=\lim_{\rho\rightarrow0^{+}}\int_{B_{2\rho}(x_j)}\chi(\varepsilon z)(u^{+})^{2_s^{\ast}}\psi_{\rho}\,{\rm d}z+\lim_{\rho\rightarrow0^{+}}\int_{B_{2\rho}(x_j)}\sum_{j\in J}\nu_j\delta_{x_j}\chi(\varepsilon z)\psi_{\rho}\,{\rm d}z\\
&=\chi(\varepsilon x_j)\nu_j,
\end{align*}
where we have used \eqref{equ7-9}. Since for any $u,v\in H^s(\mathbb{R}^3)$, there holds
\begin{equation*}
\int_{\mathbb{R}^3}\int_{\mathbb{R}^3}\frac{(u(x)-u(y))(v(x)-v(y))}{|x-y|^{3+2s}}\,{\rm d}x\,{\rm d}y=\int_{\mathbb{R}^3}(-\Delta)^{\frac{s}{2}}u(-\Delta)^{\frac{s}{2}}v\,{\rm d}x,
\end{equation*}
then it is easy to check that
\begin{equation*}
\int_{\mathbb{R}^3}(-\Delta)^{\frac{s}{2}}u_n(-\Delta)^{\frac{s}{2}}(u_n^{+}\psi_{\rho})\,{\rm d}z\geq\int_{\mathbb{R}^3}(-\Delta)^{\frac{s}{2}}u_n^{+}(-\Delta)^{\frac{s}{2}}(u_n^{+}\psi_{\rho})\,{\rm d}z.
\end{equation*}
By nonlocal Leibniz rule, we have
\begin{align}\label{equ7-13}
\int_{\mathbb{R}^3}(-\Delta)^{\frac{s}{2}}u_n^{+}&(-\Delta)^{\frac{s}{2}}(u_n^{+}\psi_{\rho})\,{\rm d}z=\int_{\mathbb{R}^3}\psi_{\rho}(-\Delta)^{\frac{s}{2}}u_n^{+}(-\Delta)^{\frac{s}{2}}u_n^{+}\,{\rm d}z\nonumber\\
&+\int_{\mathbb{R}^3}u_n^{+}(-\Delta)^{\frac{s}{2}}u_n^{+}(-\Delta)^{\frac{s}{2}}\psi_{\rho}\,{\rm d}z+\int_{\mathbb{R}^3}(-\Delta)^{\frac{s}{2}}u_n^{+} B(u_n^{+},\psi_{\rho})\,{\rm d}z.
\end{align}
Using \eqref{equ2-1} and \eqref{equ2-2}, it is easy to verify that
\begin{equation}\label{equ7-14}
\lim_{\rho\rightarrow0^{+}}\limsup_{n\rightarrow\infty}\int_{\mathbb{R}^3}u_n^{+}(-\Delta)^{\frac{s}{2}}u_n^{+}(-\Delta)^{\frac{s}{2}}\psi_{\rho}\,{\rm d}z=0
\end{equation}
and
\begin{equation}\label{equ7-15}
\lim_{\rho\rightarrow0^{+}}\limsup_{n\rightarrow\infty}\int_{\mathbb{R}^3}(-\Delta)^{\frac{s}{2}}u_n^{+} B(u_n^{+},\psi_{\rho})\,{\rm d}z=0.
\end{equation}

Therefore, by \eqref{equ7-12}--\eqref{equ7-15}, we have that
\begin{align*}
\chi(\varepsilon x_j)\nu_j&=\lim_{\rho\rightarrow0^{+}}\limsup_{n\rightarrow\infty}\int_{\mathbb{R}^3}g(\varepsilon z,u_n)u_n\psi_{\rho}\,{\rm d}z\geq\lim_{\rho\rightarrow0^{+}}\limsup_{n\rightarrow\infty}\int_{\mathbb{R}^3}\psi_{\rho}|(-\Delta)^{\frac{s}{2}}u_n^{+}|^2\,{\rm d}z\\
&=\lim_{\rho\rightarrow0^{+}}\int_{B_{2\rho}(x_j)}\psi_{\rho}\,{\rm d}\mu=\mu(\{x_j\}).
\end{align*}
Combining with \eqref{equ7-11}, we have
\begin{equation}\label{equ7-16}
\nu_j\geq\Big(\frac{\mathcal{S}_s}{\chi(\varepsilon x_j)}\Big)^{\frac{3}{2s}}.
\end{equation}

Considering $\eta_R(z):=\eta(\frac{z}{R})$ for $R>0$, where $\eta$ is a smooth cut-off function such that $\eta=0$ on $B_1(0)$, $\eta=1$ on $\mathbb{R}^3\backslash B_2(0)$, $0\leq\eta\leq1$ and $|\nabla\eta|\leq C$. Suppose that $R$ is chosen in such a way that $\Lambda'/\varepsilon\subset B_{R}(0)$. By $(g_3)$, we have that
\begin{align}\label{equ7-17}
\int_{\mathbb{R}^3}(-\Delta)^{\frac{s}{2}}u_n(-\Delta)^{\frac{s}{2}}(u_n^{+}\eta_R)\,{\rm d}z&+\int_{\mathbb{R}^3}V(\varepsilon z)(u_n^{+})^2\eta_R\,{\rm d}z+\int_{\mathbb{R}^3}\phi_{u_n}^t(u_n^{+})^2\eta_R\,{\rm d}z\nonumber\\
&=\int_{\mathbb{R}^3}g(\varepsilon z,u_n)u_n^{+}\eta_R\,{\rm d}z+o(1)\nonumber\\
&=\int_{\mathbb{R}^3\backslash(\Lambda'/\varepsilon)}\tilde{f}(u_n)u_n^{+}\eta_R\,{\rm d}z+o(1)\nonumber\\
&\leq\frac{\alpha_0}{k}\int_{\mathbb{R}^3\backslash(\Lambda'/\varepsilon)}(u_n^{+})^2\eta_R\,{\rm d}z+o(1).
\end{align}
Since
\begin{equation*}
\lim_{R\rightarrow+\infty}\limsup_{n\rightarrow\infty}\int_{\mathbb{R}^3}(u_n^{+})^2\eta_R\,{\rm d}z=\lim_{R\rightarrow+\infty}\int_{|z|>R}(u^{+})^2\eta_R\,{\rm d}z=0
\end{equation*}
and using the similar argument to \eqref{equ2-5} and \eqref{equ2-6} to deduce that
\begin{align*}
\lim_{R\rightarrow+\infty}\limsup_{n\rightarrow\infty}\int_{\mathbb{R}^3}(-\Delta)^{\frac{s}{2}}u_n(-\Delta)^{\frac{s}{2}}(u_n^{+}\eta_R)\,{\rm d}z&\geq\lim_{R\rightarrow+\infty}\limsup_{n\rightarrow\infty}\int_{\mathbb{R}^3}|(-\Delta)^{\frac{s}{2}}u_n^{+}|^2\eta_R\,{\rm d}z\\
&=\mu_{\infty},
\end{align*}
thus, by \eqref{equ7-17} and \eqref{equ7-11}, we conclude that $\mu_{\infty}=\nu_{\infty}=0$.

On the other hand, by $(g_2)$, \eqref{equ7-11} and $\mu_{\infty}=\nu_{\infty}=0$, we have that
\begin{align*}
&c+o(1)\\
&=J_{\varepsilon}(u_n)-\frac{1}{4}\langle J_{\varepsilon}'(u_n),u_n\rangle\\
&\geq\frac{1}{4}\int_{\mathbb{R}^3}|(-\Delta)^{\frac{s}{2}}u_n^{+}|^2\,{\rm d}z+\frac{1}{4}\int_{\mathbb{R}^3}V(\varepsilon z)(u_n^{+})^2\,{\rm d}z+\int_{\mathbb{R}^3}(\frac{1}{4}g(\varepsilon z,u_n)u_n-G(\varepsilon z,u_n))\,{\rm d}z\\
&=\frac{1}{4}\int_{\mathbb{R}^3}|(-\Delta)^{\frac{s}{2}}u_n^{+}|^2\,{\rm d}z+\frac{1}{4}\int_{\mathbb{R}^3}V(\varepsilon z)(u_n^{+})^2\,{\rm d}z+\int_{\mathbb{R}^3}\chi(\varepsilon z)(\frac{1}{4}f(u_n)u_n-F(u_n))\,{\rm d}z\\
&+\frac{4s-3}{12}\int_{\mathbb{R}^3}\chi(\varepsilon z)(u_n^{+})^{2_s^{\ast}}\,{\rm d}z+\int_{\mathbb{R}^3}(1-\chi(\varepsilon z))(\frac{1}{4}\tilde{f}(u_n)u_n-\tilde{F}(u_n))\,{\rm d}z\\
&\geq\frac{1}{4}\int_{\mathbb{R}^3}|(-\Delta)^{\frac{s}{2}}u_n^{+}|^2\,{\rm d}z+\frac{4s-3}{12}\int_{\mathbb{R}^3}\chi(\varepsilon z)(u_n^{+})^{2_s^{\ast}}\,{\rm d}z\\
&\geq\frac{1}{4}\mu(\{x_j\})+\frac{4s-3}{12}\chi(\varepsilon x_j)\nu_j+o(1).
\end{align*}
Thus, using \eqref{equ7-11} and \eqref{equ7-16}, one has
\begin{align*}
c\geq\frac{1}{4}\mu(\{x_j\})+\frac{4s-3}{12}\chi(\varepsilon x_j)\nu_j&\geq\frac{1}{4}\mathcal{S}_s\nu_j^{\frac{2}{2_s^{\ast}}}+\frac{4s-3}{12}\chi(\varepsilon x_j)\nu_j\\
&\geq\frac{1}{4}\mathcal{S}_s\Big(\frac{\mathcal{S}_s}{\chi(\varepsilon x_j)}\Big)^{\frac{3}{2s}}+\frac{4s-3}{12}\frac{\mathcal{S}_s^{\frac{3}{2s}}}{\chi(\varepsilon x_j)^{\frac{3-2s}{2s}}}\\
&=\frac{s}{3}\frac{\mathcal{S}_s^{\frac{3}{2s}}}{\chi(\varepsilon x_j)^{\frac{3-2s}{2s}}}\geq\frac{s}{3}\mathcal{S}_s^{\frac{3}{2s}}
\end{align*}
which leads to a contradiction. Hence \eqref{equ7-8} holds, then the claim \eqref{equ7-4} is true. Combining $\langle J_{\varepsilon}'(u_n),u_n\rangle=0$ with \eqref{equ7-2} and \eqref{equ7-4}, by standard argument, we can get $u_n\rightarrow u$ in $H_{\varepsilon}$.
\end{proof}

{\bf Acknowledgements.}
The work is supported by NSFC grant 11501403.

\end{document}